\documentclass[11pt]{article}

\usepackage{amsthm,amsfonts}
\usepackage{amsmath,amssymb, enumerate, fullpage}
\usepackage{graphicx,color,xcolor}
\usepackage{hyperref}
\theoremstyle{definition}

\DeclareMathOperator*{\argmin}{arg\,min}

\newcommand{\bc}{\mathbf{c}}
\newcommand{\bb}{\mathbf{b}}
\newcommand{\bx}{\mathbf{x}}
\newcommand{\by}{\mathbf{y}}
\newcommand{\cL}{\mathcal{L}}
\newcommand{\wT}{\widetilde{T}}
\newcommand{\resid}{\text{resid}}

\begin{document}

\title{IDENT Review: Recent Advances in Identification of Differential Equations from Noisy Data} 
\author{Roy Y. He\thanks{Department of Mathematics, City University of Hong Kong. Email: royhe2@cityu.edu.hk.},  
~Hao Liu\thanks{
Department of Mathematics, Hong Kong Baptist University. Email: haoliu@hkbu.edu.hk.},  
~Wenjing Liao\thanks{School of Mathematics, Georgia Institute of Technology, Atlanta, GA 30332-0160. Email: wliao60@gatech.edu.} 
~and Sung Ha Kang\thanks{School of Mathematics, Georgia Institute of Technology, Atlanta, GA 30332-0160. Email: kang@math.gatech.edu.}}

\date{}

\maketitle

\begin{abstract}
Differential equations and numerical methods are extensively used to model various real-world phenomena in science and engineering.  With modern  developments, we aim to find the underlying differential equation from a single observation of time-dependent data.  If we assume that  the  differential equation is a linear combination of various linear and nonlinear differential terms, then the identification problem can be formulated as solving a linear system. The goal then reduces to finding the optimal coefficient vector that best represents the time derivative of the given data.
We review some recent works on the identification of differential equations. We find some common themes for the improved accuracy: (i) The formulation of linear system with proper denoising is important, (ii) how to utilize sparsity and model selection to find the correct coefficient support needs careful attention, and (iii) there are ways to improve the coefficient recovery.
We present an overview and analysis of these approaches about some recent developments on the topic.
\end{abstract}

%\textbf{Keywords: differential equation identification, data-driven modeling, inverse problem, numerical method for PDE}

\section{Introduction}
We review the inverse problem of finding the governing  differential equation from a single noisy  observation of the solution of this differential equation.  The  motivation is to explain real world phenomenon from the given data.  When an interesting physical phenomenon happens,  by taking a video, i.e. by collecting time-dependent data, can one explain the governing equation to help us better understand the phenomenon?  To handle such cases, we do not assume  multiple realizations of the same underlying equation, but focus on finding the best equation from one single possibly noisy observation.

This question of  understanding the underlying dynamics of data has been around for a while in science and engineering.   Some early ideas include  \cite{bellman1969new} where system identification  is considered, \cite{akaike1974new}  where statistical model identification is explored, \cite{baake1992fitting} where  ODE fitting is explored for chaotic data, \cite{muller2004parameter} where parameter estimation is explored for PDE with respect to noise and temporal and spatial resolution, \cite{bongard2007automated} where a method is explored to automatically generate symbolic equations, and \cite{schmidt2009distilling} about distilling natural laws from experimental data. 
In recent years, sparse regression is added to find a smaller number of terms in model identification, e.g., 
Sparse Identification of Nonlinear Dynamics (SINDy) \cite{brunton2016discovering,rudy2019data}, learning PDE via data discovery and sparse optimization \cite{schaeffer2017learning} and Weak SINDy ODE and PDE \cite{messenger2021weakode, messenger2021weak}, just to mention a few. 

In this paper, starting from the settings of Identification of Differential Equation using Numerical Techniques (IDENT) \cite{kang2021ident}, we present developments in this direction including Robust IDENT \cite{he2022robust}, Weak IDENT \cite{tang2023weakident}, Fourier IDENT \cite{tang2023fourier} as well as GP-IDENT \cite{he2023group,he2025group}  and some related works.  
There are common themes for better identification of differential equations:

\textbf{1. Feature system and denoising}: 
First of all, the formulation of a feature system with proper denoising is important.  There are different ways to establish the feature system, e.g. direct discretization using finite difference schemes~\cite{kang2021ident,he2022robust,he2023group} or  integral equations and weak form as in \cite{tang2023weakident}.  
As shown in \cite{kang2021ident}, denoising stabilizes identification, and an improved method is proposed in \cite{he2022robust} using finite difference forms.  
In the weak form, an idea similar to the use of a low-pass filter through convolution with test functions is explored in \cite{tang2023weakident}.

\textbf{2. Support identification}:
Secondly, many real-world physical phenomena can be expressed as a small number of feature terms, and it is often beneficial to find a simple expression than a complex one.  We use sparsity-inducing models to find candidate supports \cite{dai2009subspace}.  However, it is often difficult to identify the correct differential equation just with sparsity regularization; we thus introduce various model selection methods to refine the identification result.  These model selection methods are mathematically based on the fundamentals of numerical PDE methods and their convergence, and the effects of least-square minimization.  We present time evolution error, multishooting time evolution error, cross validation as well as reduction in the residual for the purpose of model selection. 

\textbf{3. Coefficient value recovery and further refinements:} Typically, constant coefficient values are computed by using least-square fit after the correct support (i.e. the form of the equation) is found.  There are ways to improve the coefficient value recovery when using least-squares, such as, using high dynamic region and narrow fit, which helps to improve coefficient reconstruction accuracy.
Moreover, this set-up can be used to find differential equations with  coefficients that may change in time and/or space.  We present the details for methods in this direction, using group-LASSO, group Subspace Pursuit, and a patch based approach.  We also present the column normalization method which helps to stabilize the computation.

We first present the set-up of the problem in Section \ref{sec:linear}. We give details on how to deal with noise in the feature system in Section \ref{sec:denoise} with error analysis in Subsection \ref{ssec:error}.  
The main part of support recovery via sparsity regularization and model selection is presented in Section \ref{sec:support}.  We present additional refinement methods in Section \ref{sec:refinement}, and collective results for constant coefficient differential equation identification in Section \ref{sec:numerical}. 
We discuss identifying differential equation with varying coefficients in Section \ref{sec:varying}, and we conclude the paper with summarizing remarks in Section \ref{sec:summary}.

\section{Set-up: construction of the feature system}\label{sec:linear}

We let the given set of time-dependent data be 
\[ 
\{U_{i}^n | i=1, \dots, N_x \text{ and } n= 1, \dots, N_t\},
\]
where the indices $i$ and $n$ are associated with the grid point $(x_i,t_n)$ in the spatial and time domain, respectively.  We present the details in one spacial dimension, but the formulation can be generalized to high spacial dimensions.  Since there may be  noise in the observation,  we assume \begin{equation}
U_{i}^n = u(x_i, t_n) + \epsilon_i^n 
\label{eq:givedata}
\end{equation}
where $u(x_i, t_n)$ is the true data.  
To have a starting point, we assume that the governing differential equation of $u (x,t)$  is a linear combination of various linear and nonlinear differential terms, such that 
\begin{align}\label{E:general_constant}
 u_t = \mathcal{F}(x,u,u_x,u_{xx},\dots) 
= \sum_{k=1}^{N_f} c_kf_k 
\end{align}
where $f_k$'s are monomials of $u$, e.g., 1, $u$, $u_x$, $u^2$, etc, which we refer to as \textbf{feature terms}.  There are $N_f$ number of feature terms, and $c_k$'s are scalars.  The right hand of~\eqref{E:general_constant} is typically called the governing function of the data which is assumed to be unknown.

Using the given discrete time and spacial  data $\{ U^n_i\}$, the model \eqref{E:general_constant} can be approximated by a discrete linear system:
\begin{equation}\label{E:linearsystem}
 \bb = F \textbf{c} 
\end{equation}
where $\bb$ is a discrete approximation to $u_t$, and 
each column of the \textbf{feature matrix $F$} is an approximation to feature term
{\footnotesize
\begin{equation} \label{E:constant_F}
 F \approx \left( \begin{array}{ccccccc}
\vdots & \vdots & \vdots & \vdots &\vdots & \vdots &\vdots \\% & \vdots &\vdots & \vdots \\
1 & u(x_i,t_n) & u^2(x_i,t_n) & u_x(x_i,t_n) & u_x^2(x_i,t_n) & uu_x(x_i,t_n) & \cdots \\%& u_{xx}(x_i,t_n) & u_{xx}^2(x_i,t_n) & u u_{xx}(x_i,t_n) & u_x u_{xx}(x_i,t_n) \\
\vdots & \vdots & \vdots & \vdots &\vdots & \vdots &\vdots \\%& \vdots &\vdots & \vdots \\
\end{array}  \right).
\end{equation}}
This feature matrix is multiplied by \textbf{a coefficient vector}, $\textbf{c} = (c_1, c_2, \dots, c_{N_f})^T$. 
In this form, the identification of differential equation becomes identifying the vector $\textbf{c}$.  The index of non-zero coefficients, we refer to as the  \textbf{support}, gives the form of the equation, and the \textbf{coefficient value} gives the coefficients of the equation \eqref{E:general_constant}.  The goal is to find the best coefficient vector $\textbf{c}$, the support as well as the values, which best expresses the dynamics of the given data.

For the derivative approximations in the feature matrix $F$ in \eqref{E:constant_F}, in IDENT \cite{kang2021ident}, Robust IDENT \cite{he2022robust} and GP-IDENT\cite{he2023group}, finite difference method, e.g. the ENO scheme \cite{harten1997uniformly}, is used to approximate spacial derivatives with high order accuracy even if there are  discontinuities in the given data.  With noise in the observed data, denoising helps to stabilize the identification (see Sections \ref{ssec:lsma} and \ref{ssec:sdd}). 

Using an integral form, in particular a weak form \cite{gurevich2019robust,messenger2021weak,messenger2021weakode, tang2023weakident}, one obtains the following equation against the test function $\phi_{h(x_i,t_n)}$:  
\begin{equation}
      \int_{\Omega_{h(x_i, t_n)}} \phi_{h(x_i, t_n)}(x,t) \frac{\displaystyle \partial u(x,t)}{\displaystyle \partial t}dx dt = \sum_{k=1}^K c_k  \int_{\Omega_{h(x_i, t_n)}} \phi_{h(x_i, t_n)}(x,t)f_k dxdt.
      \label{eq:weakformulation1}
\end{equation}
Here $\Omega_h(x_i,t_n)$ is the support of $\phi_h(x_i,t_n)$, which is an open neighborhood of $(x_i,t_n)$ with size adjustable via the parameter $h$, see Section \ref{ssec:testfunction} for details.
The integral form in \eqref{eq:weakformulation1} shows benefits,  when each feature $f_k$ has the form 
\begin{equation}\label{E:Wfeature}
    f_k =\frac{\partial^{\alpha_k} } {\partial x^{\alpha_k}} u^{\beta_k},
\end{equation}
for some $\alpha_k\in\{0,1,\ldots,\bar\alpha\}$ and $\beta_k\in\{0,1,\ldots,\bar\beta\}$. 
First of all, considering integration by parts of \eqref{eq:weakformulation1}, the high order derivative can move to the smooth test function as 
\begin{equation}
    -\int_{\Omega_{h(x_i, t_n)}} u(x,t) \frac{\displaystyle \partial \phi_h(x,t)}{\displaystyle \partial t} dxdt =
    \sum_{k=1}^K c_k
    \int_{\Omega_{h(x_i, t_n)}} 
    (-1)^{\alpha_k}  u^{\beta_k}
    \frac{\partial^{\alpha_k} \phi_h}{\partial x^{\alpha_k}} dxdt,
    \label{e: integral form system}
\end{equation}
as long as $\phi_h$ and its derivatives up to order $\bar{\alpha}$ vanish on the boundary of $\Omega_{h(x_i,t_n)}$.  This can give higher order accuracy  compared to using finite difference methods, since the test function can be very smooth.  Secondly, with a well-designed test function, better denoising can be achieved (see Section \ref{ssec:testfunction}). The weak form in \eqref{e: integral form system} can be formulated as a linear system 
\begin{equation}
    W\bc = \bb, 
\label{e: Wc=b}
\end{equation}
where the rows are indexed by $h$.
Furthermore,  Fourier modes can be used to represent features in the frequency domain as in Fourier IDENT \cite{tang2023fourier}, and differential equations with varying coefficients can also be represented as a linear system (see Section \ref{sec:varying}).

In general, just by assuming the underlying  governing  differential equation is a linear combination of linear and nonlinear terms, this inverse problem can be represented as a linear system.  Then the identification of differential equations becomes finding the support and values of the coefficient vector $\bc$.  
While this linear system seems relatively easy to solve, the noise added on the observation  increases the complexity of the problem and yields instability.   We review denoising methods for this linear system to stabilize this process, and present related error analysis in the subsequent section.

 %%%%%%%%%%%%%%%%%%%%%%%%%%%%%%%%%%%%%%%%
\section{Denoising identification linear system}%Feature Estimation from Noisy Data}
\label{sec:denoise}  

In differential equation identification, denoising is important in setting up the feature system. For example, Table \ref{table:burgers} shows an example of Burgers’ equation $u_t = - u u_x$ with Dirichlet boundary conditions considered in \cite[Equation (15)]{kang2021ident} .  This table is showing the top two equations with the smallest errors, in particular time evolution error.  
When the given data has no noise,  identification results can be good.  When noise is present, the identified equation can be quite different from the true equation shown in the second column.   With a good denoising (in this case, least-square moving average is used), correct equation can be identified as shown in the last row.  The idea of this table is presented in \cite{kang2021ident}.

\subsubsection{Denoising the given data}
\label{ssec:lsma}

Let  $\{U_i\}$ be the given data on a one-dimensional uniform grid $\{x_i\}$ and define its five-point moving average as  $\tilde{U}_i = \frac{1}{5}\sum_{l=0,\pm1, \pm2} U_{i+l}$ for all $i$. 
At each grid point $x_i$, we determine a quadratic polynomial $p(x)= a_0 + a_1 (x-x_i) + a_2 (x-x_i)^2$ fitting the local data, which preserves the order of smoothness, up to the degree of polynomial. 
There are a few possible choices for denoising, such as (i) Least-Squares Fitting (LS): find $a_0, a_1$ and $a_2$ to minimize the functional
$ F( a_0, a_1, a_2 )=\sum_{{\rm some} \; j\; {\rm near} \; i} (p(x_j)- U_j)^2$, and 
(ii) Moving-Average Fitting (MA): find $a_0, a_1$ and $a_2$, such that the local average of the fitted polynomial matches with the local average of the data, $ {1}/{5} \sum_{l=0,\pm1, \pm2} p(x_{j+l}) = \tilde{U}_j,$ for $j=i, i\pm 1$.
In IDENT \cite{kang2021ident}, (iii) Least-Squares Moving Average (LSMA) is used: find $a_0, a_1$ and $a_2$ to minimize the functional
\[G( a_0, a_1, a_2 )=  \sum_{j=i,i\pm1, i\pm2} \left\{\left[\frac{1}{5} \sum_{l=0,\pm1, \pm2} p(x_{j+l})\right] - \tilde{U}_j\right\}^2.\]
This LSMA denoising method preserves the approximation order of data and can easily be incorporated into numerical PDE techniques. The quadratic polynomials computed by the methods above are locally third-order approximation to the underlying function. Theorem 2 in \cite{kang2021ident} states that if the given data are sampled from a third-order approximation to a smooth function, then LSMA will keep the same order of the approximation. This theorem can be easily generalized to any higher order.  
In \cite{kang2021ident}, the authors show how LSMA stabilizes the process as in Table \ref{table:burgers}, and also gives better results compared to no denoising, or using LS of MA. 

\begin{table}
\centering
\begin{tabular}{c|c|c}
\hline
 & Identified equation & error (TEE)\\
 \hline
clean data 
 & $u_t = -0.99 u u_x $ & 0.48 \\
 & $u_t = 0.1 u -0.99 u u_x$ & 1.40 \\
% & $u_t = 0.27 u$ & 78.76 \\
\hline
$8\%$ noise
& $u_t = -0.59 u^2  $ & 94.03\\
without  denoising & $u_t = 0.18 - 0.83 u^2  $ & 94.04\\
%& $u_t = -0.27 $ & 94.20\\
\hline
$8\%$ noise
& $u_t = -0.92 u u_x  $ & 25.5409\\
with denoising & $u_t = -0.22 - 0.92 u u_x $ & 29.5412\\
%& $u_t = -0.22 + 0.07 u -0.92 u u_x $ & 29.80\\
\hline
\end{tabular}
\caption{Importance of denoising: The true equation is $u_t = - u u_x$, the  equation (15) in \cite{kang2021ident}.   This table is showing the top two equations with the smallest errors (TEE) (see Section~\ref{sec:support}). When the given data are clean, the identification results can be quite good.  When noise is present, a good denoising method is needed to stabilize the process. }
\label{table:burgers}
\end{table}

In asymptotic analysis \cite{he2022asymptotic}, the authors proved that $\ell^1$-based identification method with kernel smoothing is guaranteed to asymptotically  recover the signed-support of the coefficients for   
\begin{align*} 
\mathbf{c}(\lambda)=\arg\min_{\mathbf{c}\in \mathbb{R}^{N_f}}
    \bigg\{\frac{1}{2N_x N_t}\left\|\mathbf{b}-{F}\mathbf{c}\right\|_{2}^{2}+\lambda_{N}\|\mathbf{c}\|_1\bigg\}\;.
\end{align*}
The authors employ a local quadratic regression to estimate $u_t(x_i,\cdot)$ for each fixed space point $x_i$ and use a local-polynomial with degree $p+1$ to estimate $\partial_x^pu(\cdot,t_n)$ at each temporal point $t_n$, for each degree $p=0,1,\dots,P_{\max}$. 

\subsubsection{Denoising higher order derivatives}\label{ssec:sdd}

Denoising the given data helps with identification, yet, the feature terms in $F$ have derivative terms.   Using finite difference approximation, the noise can get emphasized as derivatives are approximated.   In \cite{he2022robust}, Successively Denoised Differentiation (SDD) is proposed to better handle higher order derivatives.  
Let $S$ be a smoothing operator, e.g. Moving Least Square (MLS)~\cite{lancaster1981surfaces}, 
a weighted least squares problem is solved at each  spatial location $x_{i}$ 
\begin{align*}
&S_{(x)} \left[U\right]_i^n = p^n_i(x_i),\nonumber\\
&\text{ with}\;\; p^n_i=\argmin_{p\in P_2}\sum_{0\leq j\leq N_x}(p(x_j)-U_{j}^n)^2	\exp\left(-\frac{\|x_i-x_j\|^2}{h^2}\right),
\end{align*}
here $h>0$ is a width parameter of the kernel. The main idea of SDD is to smooth the data at each step (before and after) of the numerical differentiation. For example,
\[ \partial_{x}^{k} u 
\approx (S_{(x)}D_{x})^{k}S_{(x)}[U] \]
that spatial denoising is applied after every numerical spatial differentiation $D_x$.   This procedure can effectively stabilize the numerical differentiation.  
Theorem 2.1 in \cite{he2022robust} shows that under appropriate assumptions, for clean data, the estimated partial derivative $S_{(x)}D_x[U]$ has the same accuracy as the estimated derivative without SDD. 

Differentiation via SDD is versatile as it can be applied to any form of feature terms, such as $u u_{xx}^2$, $u_x^2 u_{xx}^3$, etc.  This contrasts with the use of integral forms where in order to take advantage of integration by parts, the feature terms should have a particular form as in \eqref{E:Wfeature}.  As shown in \cite{he2022robust}, SDD improves the feature term approximation and helps with better identification when using finite difference for the approximation of feature terms.

\subsubsection{Denoising in weak form}\label{ssec:testfunction}

In the weak form~\eqref{e: integral form system}, the integration with test functions gives rise to a denoising effect.  We find that the following test function is a good option \cite{messenger2021weak}:
\[
    \phi(x,t) =  \left(1 - \left(\frac{x}{m_{x}\Delta x}\right)^2\right)^{p_x}\left(1-\left(\frac{t}{m_t\Delta t} \right)^2\right)^{p_t}, \;\; \; 
    |x|\le m_x\Delta x, \ |t|\le m_t\Delta t
\]
for $\phi_h$ as a translation of $\phi$  for  $i = 1,..,N_x, n = 1,...,N_t$.  Here $p_x$ and $p_t$ denote the polynomial order which is related with the smoothness of $\phi$ in terms of $x$ and  $t$. The test function satisfies  $\int_{\Omega_{h(x_i,t_n)}} \phi_h(x,t)dx dt  = 1$ and $ \phi_h(x,t) = 0 \; \text{for} \; (x,t) \in \partial \Omega_{h(x_i, t_n)} $, with $\phi(x,t)$ localized around $(x_i,t_n)$ and is supported on $\Omega_{h(x_i,t_n)}  = [x_{i} - m_x \Delta x, x_{i} + m_x \Delta x ] \times [t_{n} - m_t \Delta t, t_{n} + m_t \Delta t]$ for some positive integers $m_x$ and $m_t$.  
The weak features %$ w_{h(x_i,t_n)}$ 
in \eqref{e: Wc=b} can be written in a convolution form $U * \frac{\partial^{\alpha_l} }{\partial x^{\alpha_l}}\phi$ and calculated through Fast Fourier Transform in terms of $\mathcal{F}^{-1}\left(\mathcal{F}(U) \circ \mathcal{F} \left( \frac{\partial^{\alpha_l} }{\partial x^{\alpha_l}}\phi \right) \right)$, where $\circ$ denotes point-wise multiplication, and  $p_x$, $p_t$, $m_x$ and $m_t$ are carefully chosen to give a denoising effect depending on the frequency of the given data as detailed in \cite{messenger2021weak}.

In Fourier IDENT \cite{tang2023fourier}, denoising is  applied in frequency domain for each feature. The denoising effect is analogous to applying Gaussian kernel convolution to noisy data in the physical domain, but in Fourier IDENT,  it is performed in the frequency domain.

\subsection{Error analysis}\label{ssec:error}

We summarize the error analysis for identification of differential equations.   When  finite difference is used to approximate derivatives, the identification is harder when the underlying differential equation is of a high order, while using the weak form reduces this effect.

In IDENT \cite{kang2021ident}, the differential features  in $F$ and $\bb$ in  \eqref{E:linearsystem} are empirically approximated by finite difference schemes.  Let {$F$ and $\bb$ be the empirical feature matrix and time derivative vector in \eqref{E:linearsystem} approximating the exact feature matrix $F^*$ and time derivative vector   $\bb^*$ respectively. Denote the underlying true coefficient vector by $\bc^*$. Then the empirical linear system ${F} \bc^* \approx \bb $
has error
$$ \mathbf{e} =  \bb^* -{\bb} + ({F}-F^*)\mathbf{c}^*,  $$ where 
$\bb^* -{\bb}$ is resulted from the approximation of $u_t$, and  $({F}-F^*)\mathbf{c}^*$ is resulted from the approximation of the spatial derivatives of $u$. It was shown in~\cite{kang2021ident} that 
\begin{equation}
\label{E:enoise}
\|\mathbf{e}\|_{L^2}  \le O\left(\Delta t + \Delta x^{p+1-r}+ \underbrace{\frac{\delta t + \delta x^q}{\Delta t}+ \frac{\delta t + \delta x^q }{\Delta x^r}}_{\text{errors from data generations}} + \underbrace{\frac{\sigma}{\Delta t} + \frac{\sigma}{\Delta x^r}}_{\text{measurement noise}} \right), 
\end{equation}
where $r$ is the highest order of the spacial derivatives in \eqref{E:general_constant}.
This considers the numerical simulation error and the noise. 
Suppose the clean data are numerically simulated for the PDE in \eqref{E:general_constant} by a $q$th-order method with time step $\delta t$ and spatial spacing $\delta x$, and the measurement noise  is independently drawn from the normal distribution with mean $0$ and variance $\sigma^2$. Then, the noise in \eqref{eq:givedata} is
$$\epsilon _i^n = O (\delta t + \delta x^q + \sigma),$$
and we denote the size of the sampling grid in the spatial and time domain by $\Delta x$ and $\Delta t$, respectively. 
The error formula \eqref{E:enoise} suggests that finite difference is sensitive to measurement noise.    Gaussian noise with mean $0$ and variance $\sigma^2$ results in $O(\sigma/\Delta t + \sigma/\Delta x^r)$ in the error formula and higher-order PDEs are more sensitive to measurement noise than lower-order PDEs.

%%%%%%%%%%%%%%%%%%%%%%%%%%%%%%%%%%%%%

Using the Successively Denoised Differentiation (SDD) in Robust IDENT \cite{he2022robust} in a differential form also results in residual error depending on the derivatives of the governing equations. Assume SDD of order $s$ is used for the differential features in the system \eqref{E:linearsystem}. It gives rise to the following error at the query point $(x_i, t_n)$:
\begin{equation}\label{E:robustError}
\left({{F}}\bc^* - {\bb} \right)_{(x_i,t_n)}
=   O\left(\Delta t^{s} 
    + \Delta x^{s - \max_l \alpha_l + 1} + \frac{\sigma}{\Delta t}   +  \frac{\sigma}{\Delta x^{\min_l \alpha_l}}  \right) 
\end{equation}
where $\bc^*$ denotes the true coefficient, the noise  satisfies $\epsilon_i^n \sim O(\sigma)$, and $\max_{l}\alpha_l$ and $\min_{l}\alpha_l$ are the highest and lowest order of derivatives in the underlying PDE, respectively. In \cite{he2022robust}, various errors are considered in addition, such as $e(u):=U_i^n-u(x_i,t_n)$ and  $e(u_t):=D_tU_i^n-u_t(x_i,t_n)$ which is decomposed into response error, regression error, coefficient error and system error.  

For Weak IDENT \cite{tang2023weakident}, if we denote the support of the true coefficient vector  by ${\rm Supp}^* = \{ k: \ c_k^*\neq 0\}$, the  associated integral formulation in \eqref{e: integral form system} with the true coefficients from the true support becomes
\[
    \int_{\Omega_h}u(x,t) \frac{\displaystyle \partial \phi_h(x,t)}{\displaystyle \partial t} dxdt +
    \sum_{l\in \text{Supp}^*} (-1)^{\alpha_l} c_l^* \int_{\Omega_h}
     u^{\beta_l}(x,t)
    \frac{\partial^{\alpha_l} \phi_h(x,t)}{\partial x^{\alpha_l}}  dxdt = 0.
\]
Then the error for the discretized system in \eqref{e: Wc=b} becomes 
\begin{align}
    \mathbf{e}  \;\;\;  = \;\;\;  {{W}}\bc^* -\bb \;\;\;  = \;\;\;  \mathbf{e}^{\rm int} + \mathbf{e}^{\rm noise} 
    \label{e: error W_hat c - h_hat}
\end{align}
where we decompose the error to $\mathbf{e}^{\rm int}$ the numerical integration error of the noise free data and $\mathbf{e}^{\rm noise}$ the error due to noise.  It has been shown in \cite{messenger2021weak} that $\mathbf{e}^{\rm int} = \mathcal{O}((\Delta x \Delta t)^{q+1})$, where $q$ is the order of the numerical  integration.   An estimate of the error $\mathbf{e}^{\rm noise}$ is provided in   \cite[Theorem 1]{tang2023weakident} to be 
\begin{equation*}
    \|\mathbf{e}^{\rm noise}\|_\infty \le \left(\bar{S}^*   \max_h|\Omega_h|\right)   \epsilon + O\left(\epsilon^2\right)
\end{equation*}
where $\Omega_h$ is the support of the test function $\phi_h$, $|\Omega_h|$ denotes the area of $\Omega_h$, and  
\[\bar{S}^* =\max_h \sup_{(x_j,t_k) \in \Omega_h}
\bigg|
\sum_{l\in \text{Supp}^*}(-1)^{\alpha_l}{c_l}^*
{\beta_l}(U_j^k)^{ \beta_l - 1}
\frac{\partial^{\alpha_l} \phi }{\partial x^{\alpha_l}}(x_j, t_k)
-  \frac{\displaystyle \partial \phi}{\displaystyle \partial t}  (x_j, t_k) 
\bigg|.\]
The total error $
\mathbf{e}   $ in (\ref{e: error W_hat c - h_hat}) is bounded by 
\begin{equation}\label{eq:errorBound}
\|\mathbf{e}\|_\infty \leq O ((\Delta x \Delta t)^{q+1}) + \bar{S}^* |\Omega_h| \epsilon + O(\epsilon^2) \end{equation}
where $\epsilon = O(\sigma)$ represents the noise level. The error in \eqref{eq:errorBound} depends on the accuracy of the numerical integration and the size of the integration region, but is otherwise independent of the differential error, i.e., the error of the weak form in  (\ref{eq:errorBound}) does not depend on the order of the derivative of the governing equation compared to the error terms in \eqref{E:enoise} and \eqref{E:robustError}.  We observe that the error for the discretized linear system in the weak form is significantly smaller than the error using finite difference methods for the differential forms, when the feature terms are of form \eqref{E:Wfeature}.

\subsection{Theoretical results on identifiability}

% We present analysis on when one can identify the PDE from the given data, that is, what ensures that the observed trajectory uniquely determines a governing differential equation. 
We present theoretical results on conditions when differential equations can be uniquely identified  from single trajectories. Generally speaking, when the observed trajectory data contains richer dynamic information, the underlying PDE is easier to identify.

In~\cite{he2024much}, it was shown that a linear PDE with constant coefficients can be uniquely identified from observations at two time points if and only if the solution contains enough Fourier modes. This  confirms empirical evidence that identification is more robust when the observed trajectories exhibit sufficient variation. 
% one can transform it into the Fourier domain and show that the underlying differential operator $\cL$ can be identified by one trajectory at two different instants
Define the Fourier transform of $u$ with respect to the space variable 
$$\widehat{u}(\zeta, t) = (2\pi)^{-d/2}\int_{\Omega} e^{-i\zeta\cdot x} u (x, t)\,dx.$$
The PDE $\partial_t u = \cL u$ is converted to an ODE for each frequency $\zeta$,
\[
	\partial_t\widehat{u}(\zeta, t) = -(2\pi)^{-d/2}\sum_{|\alpha|=0 }^np_\alpha(i\zeta)^{\alpha}\widehat{u}(\zeta, t), \] 
whose solution is
\begin{align*}
	\widehat{u}(\zeta, t) =\widehat{u}(\zeta,0) \exp\left(-(2\pi)^{-d/2}\sum_{|\alpha|=0}^np_\alpha (i\zeta)^\alpha t\right)\;.
\end{align*}
Suppose there is a $\zeta\in\mathbb{R}^d$ such that $\widehat{u}(\zeta,0) \neq 0$, then for any $t_2> t_1>0$, we have
{\small \begin{align}\nonumber
	\frac{\widehat{u}(\zeta, t_2)}{\widehat{u}(\zeta, t_1)} =\exp\left(-(2\pi)^{-d/2}\sum_{|\alpha|~\text{even}}^n p_\alpha (i\zeta)^\alpha (t_2-t_1)\right) \exp\left(-(2\pi)^{-d/2}\sum_{|\alpha|~\text{odd}}^n p_\alpha (i\zeta)^\alpha (t_2-t_1)\right).
\end{align}}
By denoting $c_\alpha = p_\alpha i^{|\alpha|}$ when $|\alpha|$ is even, and $c_\alpha = p_\alpha i^{|\alpha|-1}$ when $|\alpha|$ is odd, we associate the  $(\zeta, t)$ pair with the following decoupled system
\begin{align}
		\frac{(2\pi)^{d/2}}{t_2-t_1}\log\left(\left|\frac{\widehat{u}(\zeta, t_2)}{\widehat{u}(\zeta, t_1)}\right|\right)&=  -\sum_{|\alpha|\leq n,~|\alpha|~\text{even}}c_\alpha \zeta^\alpha\label{eq_even},\\
		\frac{(2\pi)^{d/2}}{t_2-t_1}\text{Arg}\left( \frac{\widehat{u}(\zeta, t_2)}{\widehat{u}(\zeta,t_1)}\right) &= -\sum_{|\alpha|\leq n,~|\alpha|~\text{odd}}c_\alpha \zeta^\alpha\label{eq_odd}\,.
\end{align}
This shows that there exists a unique set of parameters $p_\alpha$ such that $\partial_t u = -\cL u$ if and only if~\eqref{eq_even} and~\eqref{eq_odd} admit unique solutions for $c_\alpha$, which are coefficients of two polynomials. If $t_2-t_1>0$ is small enough, the phase ambiguity in ~\eqref{eq_odd} is removed. 
Theorem 3.2 in \cite{he2024much} establishes the concrete relation between the spectrum of the solution and the identifiability of linear PDEs that it is closely related to the number of different Fourier modes contained in the observed trajectory. 

In \cite{he2024much}, the information content of a solution trajectory is characterized, and their relations with the identifiability was analyzed for two different types of operators where $\cL$ is, 1) a strongly elliptic operator, and 2) a first-order hyperbolic operator. 
This characterization is dual to the Kolmogorov $n$-width \cite{kolmogoroff1936beste} of the solution trajectory as a family of functions in $L^2(\Omega)$ parameterized by time $t$.  
When $\cL$ is the Laplace operator, the solution stays close to a low dimensional space since $u(\cdot, t)$, $t>0$ is analytic in space due to the smoothing effect.  It is shown in \cite[Theorem 2.8]{he2024much} that for a general elliptic operator $\cL$, all snapshots of any single trajectory $u(x, t)$, i.e., $u(\cdot,t)$, $0<t<T$, stays $\epsilon$ close to a linear space of dimension at most of the order $O(|\log\epsilon|^2)$. This implies the intrinsic difficulty in a direct approximation of the mapping, $u(x,t)\rightarrow u(x,t+\Delta t)$ for small $\Delta t$, i.e., the generator, for a parabolic differential operator by a single trajectory. 
On the other hand, it is proved in \cite[Lemma 2.11]{he2024much} that if $\cL$ is a first-order hyperbolic operator, the data space spanned by all snapshots of a single trajectory stays $\epsilon$ close to a linear space of dimension  $O(\epsilon^{-\gamma})$, where $\gamma$ depends on the regularity of the initial data.  Numerical examples are presented in \cite{he2024much} which validated the theory.

\section{Support identification and model selection}\label{sec:support}

One of the core problems in the identification of differential equations  is to find the coefficient vector $\bc$, and we illustrate the details in this section.

Many real world physical phenomena can be expressed as a small number of feature terms, and it is  beneficial to find a simple expression than a complex one.  Motivated by this observation, we add sparsity to find the active support of the coefficient $\textbf{c}$, i.e. finding the form of the governing equation.
Different kinds of sparsity regularization are used for differential equation identification. For example, \cite{schaeffer2017learning} utilizes the $\ell^1$ regularization and  the  Douglas–Rachford algorithm to solve it. The SINDy uses the $\ell^1$ regularization in \cite{brunton2016discovering} and the $\ell^2$ regularization in \cite{rudy2017data}, and proposes the  sequential threshold based least square and ridge regression method to solve the problems.

How one utilizes sparsity to find the correct support needs careful attention, and often using a simple sparsity promoting method may not be enough. For example, in IDENT~\cite{kang2021ident}, the $\ell^1$ regularization via LASSO~\cite{tibshirani1996regression} is used, yet, it gives many variations in support identification. Time Evolution Error (TEE) is proposed to stabilize the process.  In Robust-IDENT~\cite{he2022robust}, $\ell^0$-type constraint is used to explicitly impose the number of features in the candidate model, in particular, subspace pursuit \cite{dai2009subspace} greed algorithm is used.  In Robust-IDENT~\cite{he2022robust}, the best equation is found for all  sparsity levels, then Multishooting Time Evolution Error (MTEE) or Cross Validation is used to choose the best one among the candidate equations of different sparsity levels. This $\ell^0$ constrained framework is also used in Weak IDENT \cite{tang2023weakident} and Fourier-IDENT \cite{tang2023fourier}.   
We use $\ell^0$ or $\ell^1$ regularization as a pre-processing method, then choose the one among a few choices using model selection.   New criteria for selecting the best candidate PDE, which we refer to as model selection, are presented in subsection \ref{ssec:model}.  We note these model selection methods can be used with any sparse regularization.

\subsection{Sparse regression using $\ell^1$ and $\ell^0$}
\label{ssec:sparse}

We discuss two types of sparse regression formulation of PDE identification: the $\ell^1$ regularization used in \cite{brunton2016discovering,kang2021ident,schaeffer2017learning} and the $\ell^0$ constraint used in \cite{he2022robust,tang2023weakident,tang2023fourier}.

\subsubsection{$\ell^1$ regularization and LASSO}
\label{sss:lasso}

To enforce sparsity in the  coefficient vector $ \mathbf{c}$ for the recovered equation, $\ell^1$ regularization was explored in~\cite{kang2021ident}. For PDEs with constant coefficients, this leads to the standard Least Absolute Shrinkage and Selection Operator (LASSO) \cite{tibshirani1996regression}, which was used in \cite{brunton2016discovering,kang2021ident,schaeffer2017learning} for PDE identification:
\begin{equation}
\label{eqLASSO}
\widehat{\mathbf{c}}(\lambda) ={\textstyle \arg \min_{\bc}} \left\{ \frac{1}{2} \| \bb - F\textbf{c} \|^2_2 + \lambda \|\bc\|_{1} \right\}.
\end{equation}
Here $\lambda>0$ is a balancing parameter between the first fitting term and the second regularization term.  
The number of zeros in the solution $\widehat{\mathbf{c}}(\lambda)$ increases as $\lambda$ grows. That is, larger $\lambda$ leads to candidate models with fewer feature terms. 
Using the incoherence property proposed in \cite{donoho2001uncertainty}, and following the ideas in \cite{fuchs2004sparse,tropp2004just,tropp2006just},  a recovery theory for the model~\eqref{eqLASSO} is established in \cite{kang2021ident}. 
To measure the correlation between the $j$-th and the $l$-th column of $F$, the pairwise coherence   
\[ 
\mu_{j,l}(F) = \frac{|\langle F[j] , F[l] \rangle|}{\|F[j]\|_2 \|F[l]\|_2}\;,j,l=1,2,\dots,N_f,
\]
is used, where $F[j]$ stands for the $j$th column of the feature matrix $F$, and the mutual coherence of $F$ is given as in \cite{donoho2001uncertainty}:
\[
\mu(F) = \max_{j \neq l} \mu_{j,l}(F) = \max_{j\neq l} \frac{|\langle F[j] , F[l] \rangle|}{\|F[j]\|_2 \|F[l]\|_2}.
\] 
The smaller $\mu(F)$, the less correlated are the columns of $F$, and $\mu(F) = 0$ if and only if the columns are orthogonal. By the classical results in sparse signal recovery, LASSO recovers the correct coefficients if $\mu(F)$ is sufficiently small. 
Theorem 1 in \cite{kang2021ident} shows that LASSO will give rise to the correct support when the empirical feature matrix $F$ is incoherent, i.e. $\mu(F) \ll 1$, and all underlying coefficients are sufficiently large compared to  noise. 
When the empirical feature matrix is coherent, i.e., some columns of $F$ are correlated, it has been observed that $\widehat{\mathbf{c}}(\lambda)$ are usually supported on ${\rm supp}(\mathbf{c})$ and the indices that are highly correlated with ${\rm supp}(\mathbf{c})$ \cite{fannjiang2012coherence}. 

While $\ell^1$ minimization is effective in finding a sparse vector, $\ell^1$ minimization alone is often not enough: zero coefficients in the true PDE may become non-zero in the minimizer of $\ell^1$,  if active terms are chosen by a thresholding, results are sensitive to the choice of thresholding parameter, and the balancing parameter $\lambda$ can affect the results. If some columns of the empirical feature matrix $F$ are highly correlated, LASSO is known to have a larger support than the ground truth \cite{fannjiang2012coherence}.  In \cite{kang2021ident}, Time evolution error (TEE) is proposed as a model selection technique, see Section \ref{sss:tee}.

\subsubsection{$\ell^0$ regularization and Subspace Pursuit} \label{sss:sp}

Compared to (\ref{eqLASSO}), sparsity can be controlled more closely by considering \begin{equation}
\label{robustIdent}
\widehat{\mathbf{c}}(k) ={\textstyle \arg }\min_{\bc}  \frac{1}{2} \| \bb -  F\textbf{c} \|^2_2  ,~\mbox{ subject to }\|\mathbf{c}\|_{\ell^0} = k,
\end{equation}
where $\|\bc\|_{\ell^0}=\#\{c_i\neq 0|i=1,2,\dots,N_f\}$ counts the number of non-zero elements, and $k$ is an integer with $1\leq k\leq N_f$. Note that by setting the parameter $k$, the model~\eqref{robustIdent} yields a candidate model with exactly  $k$ active features, thus, one has direct control on the model sparsity.  This shows advantages over~\eqref{eqLASSO}, which controls the sparsity level only implicitly, and it is difficult to decide what is the appropriate range for $\lambda$.  However, the problem~\eqref{robustIdent} is known to be NP-hard due to the $\ell^0$ pseudo-norm. 

In Robust IDENT \cite{he2022robust}, the authors proposed a framework based on Subspace Pursuit (SP) \cite{dai2009subspace}. 
SP is a greedy algorithm for sparse regression with given sparsity. It takes the sparsity and linear system as input, and outputs a vector of the given sparsity that best fits the linear system, allowing direct control of the sparsity of the reconstructed coefficient. Let $k$ be a positive integer and denote $\mathbf{b} = D_t U$. For a fixed sparsity level $k$, SP$(k;F,\mathbf{b})$ returns a $k$-sparse vector whose support is selected in a greedy fashion. When the given sparsity $k$ is the true sparsity and under certain conditions of the matrix $F$, such as the restricted isometry property (RIP), it is proved that SP gives the exact solution \cite{dai2009subspace}. 

Among all different sparsity level of identified equations, one best equation should be chosen. 
In~\cite{he2022robust}, Multishooting time evolution (MTEE) as well as cross-validation is used for this purposes. Cross-validation is commonly used in statistics for the choice of parameters in order to avoid overfitting~\cite{stone1978cross}, see Sections \ref{sss:mtee} and \ref{sss:cv}.

\subsubsection{Trimming feature terms}

Either by $\ell^1$ minimization or by SP, a set of candidate equations, i.e., a set of nonzero supports, are generated.  In these models, some features may give very small contributions to $u_t$, and in Weak-IDENT \cite{tang2023weakident}, trimming is  proposed to remove such features for each sparsity level $k$ after applying SP.  This reduces ambiguities in support identification. 

Consider a candidate model with sparsity $k$. Denote the coefficients of its active features by $\tilde{\bc}\in \mathbb{R}^k$ and the feature matrix corresponding to these features by $\tilde{F}$.  We define a  \textbf{contribution score} $s_i$ of each feature as 
\begin{equation*} %\label{eq:CScore}
s_i =  \frac{ n_i }{\max_{i \leq k} n_i} 
\quad \text{where} \quad 
n_i =  || \tilde{F}_i ||_2 |\bar{c}_i|, \quad  i=1,2,\dots,k , \end{equation*} 
where $\bar{c}_i$ is the $i$-th component of $\tilde{\bc}$, and $\tilde{F}_i$ is the $i$-th column of error normalized $\tilde{F}$ as in \eqref{eq:Wtilde} in Section \ref{ssec:normalization}. Since $s_i$ is normalized by the maximum value of $n_i$, $s_i$ gives the score of the contribution of the $i$-th feature relative to the contribution of the feature with the largest contribution. 
The feature term with the smallest $s_i$ which is below $\rho$, i.e., $s_i < \rho $, is trimmed, and this process is repeated by trimming features one by one until no feature can be trimmed. In \cite{tang2023weakident}, $\rho = 0.05$ is used and trimming is used for $F_{\rm narrow}$, see Section \ref{sss:narrow}. 

In Fourier-IDENT \cite{tang2023fourier}, group trimming is proposed to improve the efficiency.
Group trimming remove a set of the least relevant features for a given candidate model using a threshold $\rho$. 

\subsection{Model selection methods} \label{ssec:model}

Using $\ell^1$ or $\ell^0$ regularization is often not enough to identify the support of the coefficient, yet it gives candidate sets.  From these candidate sets, additional model selection methods are applied to find the optimal differential equation.  Time Evolution Error (TEE) is introduced in  \cite{kang2021ident}, Multi-shooting TEE (MTEE) and Cross Validation (CV) in \cite{he2022robust}, Reduction in Residual (RR) in \cite{he2023group, he2024much}, and an energy based on the residual and identification stability in \cite{tang2023fourier}.  These methods are fundamentally based on numerical PDE methods and  convergence, and least-square fit.    These model selection methods helps to find the optimal equation among different candidate equations.

\subsubsection{Time Evolution Error (TEE)} \label{sss:tee}

The main idea of Time Evolution Error (TEE) is to choose the best equation which gives the minimum dynamic error. 

For a given $\lambda$, let $\mathcal{A}$ be a set of the candidate feature index given by solving (\ref{eqLASSO}).
Consider every subset of indices $\Omega \subseteq \mathcal{A}$.  For each $\Omega= \{j_1, j_2, \ldots, j_k\}$,  find the coefficient vector  $\widehat{\mathbf{c}} = (0, 0, \widehat{c}_{j_1}, 0, \dots, \widehat{c}_{j_k}, \dots )$ by a  least-square fit such that $\widehat{\mathbf{c}}_{\Omega} = F_{\Omega}^\dagger \widehat \bb$ and $\widehat{\mathbf{c}}_{\Omega^\complement} = \mathbf{0}$, where $F_{\Omega}$ denotes the submatrix of $F$ whose column indices are in $\Omega$, and $F_{\Omega}^\dagger$ denotes the psudinverse of $F_{\Omega}$. 
Using these coefficients, the differential equation is constructed
$$
u_t=\sum_{j\in \Omega} \widehat{c}_jf_j,
$$
that is, from $\mathcal{A}$ all possible differential equations are considered.  Each candidate equation is numerically time evolved starting from the given initial data.   It is crucial to use a smaller time step $\widetilde{\Delta t}\ll \Delta t$ than the given data, where $\Delta t$ is the time spacing of the given data. 
Then the time evolution error for each $\widehat{\textbf{c}}$ is computed:
\begin{equation}\label{E:tee}
\text{ TEE} (\widehat{\textbf{c}}) := \sum_{i=1}^{N_x} \sum_{n=1}^{N_t} |\widehat{U}_i^n - U_i^n| \Delta x \Delta t,
\end{equation}
where $\widehat{U}_i^n$ is the numerically time evolved solution at $(x_i,t_n)$ of the PDE with the coefficient $\widehat{\textbf{c}}$.  The equation with the smallest TEE is chosen as the identified differential equation. 
There are two fundamental ideas behind TEE: (1) For nonlinear PDEs, it is impossible to isolate each term separately to identify each coefficient, thus, any realization of PDE must be understood as a set of terms.  (2) If the underlying dynamics are identified by the true PDE, any refinement in the discretization of the time domain should not make the numerical solution deviate from the given data.  This is the fundamental principle of consistency, stability and convergence of a numerical scheme, {and the reason for choosing a smaller time step $\widetilde{\Delta t}\ll \Delta t$.  The main effect of TEE is to evolve the numerical error from the wrongly identified differential terms.

\subsubsection{Multi-Shooting Time Evolution Error (MTEE)} \label{sss:mtee}

TEE is improved by considering multiple shooting location in Robust IDENT \cite{he2022robust}. Similar to the idea of mutishooting in solving ODE problems, Multi-Shooting Time Evolution Error (MTEE) evolves a candidate PDE from multiple time locations with a time step $\widetilde{\Delta t}\ll\Delta t$ using the forward Euler scheme for a time length of  $w\Delta t$, where $w$ is a positive integer. 
Let $\widehat{U}^{(n+w)|n}$ be the numerical solution of the candidate PDE at the time $(n+w)\Delta t$, which is evolved from the initial condition $U^n$ at time $t_n = n\Delta t$.
The MTEE is defined as
\begin{align}\label{E:mtee}
\text{MTEE}(\widehat{\mathbf{c}};w) = \frac{1}{N-w}\sum_{n=0}^{ N-1-w}\|\widehat{U}^{(n+w)|n}-U^{n+w}\|_2\;.
\end{align}
MTEE is more robust against noise in comparison with TEE in \eqref{E:tee}. If $w \ll N$,  the noise in the initial condition accumulates for a smaller amount of time in MTEE, which helps to stabilize numerical solvers.
Also, since each time evolution in the multi-shooting is independent, the computation of MTEE can be parallelized.

\subsubsection{Cross-validation Estimation Error (CEE) } \label{sss:cv}

Cross validation error is also used in Robust IDENT \cite{he2022robust} in addition to MTEE to select the best candidate PDE, and it is also used for Weak-IDENT \cite{tang2023weakident}. 

Let $\mathcal{A}$ be the set of index of a candidate differential equation. For some fixed ratio parameter $\alpha\in(0,1)$, we split the rows of $\bb$ and $F$ into two groups indexed by $\mathcal{T}_1$ and $\mathcal{T}_2$. The set $\mathcal{T}_1$ contains the first $\alpha$ portion of all rows which is used to estimate the coefficient vector by solving the following least squares problem:
\begin{align*}
\widehat{\mathbf{c}}=\argmin_{{\mathbf{c}}\in \mathbb{R}^{N_f} \text{such that}\  \mathbf{c}_{\mathcal{A}^{\complement}}=0} \|[F]_{\mathcal{A}}^{\mathcal{T}_1}\mathbf{c}_{\mathcal{A}}-[\bb]^{\mathcal{T}_1}\|_2^2\;, \end{align*}
where $[F]_{\mathcal{A}}^{\mathcal{T}_1}$ denotes the submatirx of $F$ formed by the rows indexd by $\mathcal{T}_1$ and columns indexed by $\mathcal{A}$.  The set $\mathcal{T}_2$ contains the rest of $\bb$ and $F$ which is used to validate the candidates:
\begin{align}
\mathrm{CEE}(\mathcal{A};\alpha, \mathcal{T}_1,\mathcal{T}_2)=\|[\bb]^{\mathcal{T}_2}-[F]^{\mathcal{T}_2}\widehat{\mathbf{c}}^{(k)}\|_2\;.\label{eq.CEE}
\end{align}
Here CEE stands for Cross-validation Estimation Error and it is an effective measure for consistency.  If a correct support is identified, the coefficient vector obtained from the data in $\mathcal{T}_1$ should be compatible with the data in $\mathcal{T}_2$.  Theorem 3.1 in \cite{he2022robust} shows that  CEE is guaranteed to be small provided with the correct support and sufficiently high resolution in time and space. 
Numerical experiments in~\cite{he2022robust} demonstrate that the splitting strategy of $\mathcal{T}_1$ and $\mathcal{T}_2$ does not affect the results. 
Compared to TEE \eqref{E:tee} and MTEE \eqref{E:mtee}, Cross Validation does not involve any numerical evolution of the candidate PDE, so the computation of Cross Validation is faster.

\subsubsection{Reduction in Residual (RR)}
\label{sec.RR}
In GP-IDENT \cite{he2023group}, Reduction in Residual (RR) is proposed to select the optimal candidate PDE. Suppose we have a pool of candidate PDE with distinct sparsity varying from 1 to $N_f$. Since models with a larger sparsity usually have a smaller residual, in order to balance the residual error and model sparsity, RR computes relative reduction of residual as the sparsity increases. 

For each sparsity $k$, the coefficients of the corresponding candidate model are computed by the  least square fitting, and are denoted as $\widehat{\bc}^k$. Denote the squared residual of this candidate as
$$
R_k=\|F\widehat\bc^k-\bb\|_2^2.
$$ 
Let $N_{RR}\geq 1$ be a fixed integer. For $k=1,\dots, N_f-N_{RR}$, we compute (RR) as  
 \begin{align}
 s_k = \frac{R_{k}-R_{k+N_{RR}}}{N_{RR} R_1}\;,\;\; k=1,\dots, N_f-N_{RR}.  \label{eq:RR}
 \end{align}
This measures the average reduction of residual error as the sparsity level $k$ increases. A small value in $s_k$ means there is a marginal gain in accuracy as sparsity level gets bigger than $k$.  
Here, using $N_{RR}=1$ is not reliable:  the index set $\mathcal{A}_k$ of the active features for the $k$-th candidate may not be a subset of $\mathcal{A}_{k+1}$, i.e., $R_k-R_{k+1}$ may be negative.  By using the average in~\eqref{eq:RR}, we suppress the impact of fluctuation and improve the stability of model selection. 

When the value $s_k$ is already small, we choose the smallest sparsity $k$, rather than choosing $k$ with the smallest $s_k$.  We introduce a  threshold $\rho>0$, and pick the optimal sparsity as follows:
\[
 k^*=\min\{k:1\leq k\leq N_f-N_{RR}, s_k<\rho\}. 
 \]
This is the smallest sparsity index $k$ where $s_k$ is below $\rho$. The motivation of this criterion is to find the simplest model, where RR does not significantly reduce further by considering more complex models.  For the least square fitting, as more terms are added, the error always reduces. RR helps keep simplest model being independent to increasing level of complexity with increasing sparsity level $k$.   The authors in~\cite{he2023group} find that GP-IDENT is not sensitive to the choice of $N_{RR}>1$ and $\rho$, which are chosen as $N_{RR}=5$ and $\rho =0.015$ in~\cite{he2023group}. 

RR is used in \cite{he2023group,tang2025wgidentweakgroupidentification}, and in \cite{he2024much}, a comparison among residuals of different candidate models (RRC)  is proposed. These methods are particularly useful for identifying differential equations with varying coefficient as in Section \ref{sec:varying}, as their feature matrices have larger sizes.  

\vspace{0.5cm}
In summary, these model selection methods help to  stabilize the support identification.  The main ideas are based on the  convergence of numerical PDE solvers, e.g., TEE and MTEE, and using least-squares fitting properly, e.g. CEE and RR by balancing the number of terms and the residual error.  Some of these ideas can be directly applied to all possible choice of candidate equations, but using sparsity enforcing step helps to reduce the computation. 
In Table \ref{tab_summary}, we summarize the linear system, sparsity and model selection. IDENT \cite{kang2021ident} and Robust-IDENT \cite{he2022robust} uses finite difference for feature matrix approximation with LSMA denoising and SDD, while Weak-IDENT \cite{tang2023weakident} uses weak form.  For IDENT \cite{kang2021ident}, $\ell^1$ with TEE is used for constant coefficient differential equation, for Robust-IDENT \cite{he2022robust}, SP with MTEE and CEE are used, and for Weak-IDENT \cite{tang2023weakident}, SP with Narrow fit and trimming is used.  Identification of  differential equations with varying coefficient is presented in Section \ref{sec:varying}.  

\begin{table}
\centering
\begin{tabular}{lcccc}
\hline
\textbf{Method}&\textbf{Feature system}& \textbf{Sparsity}& \textbf{Model selection}\\\hline
IDENT~\cite{kang2021ident}&Differential form&$\ell^1$ and $\ell^{2,1}$ &TEE,BEE\\\hline
Robust-IDENT~\cite{he2022robust}&Differential form&$\ell^0$-constraint& MTEE,CEE\\\hline
Weak-IDENT~\cite{tang2023weakident}&Weak  form&$\ell^0$-constraint&CEE, Trim\\\hline
Fourier-IDENT~\cite{tang2023fourier}&Fourier expansion&$\ell^0$-constraint& CEE, Group Trim\\\hline
GP-IDENT~\cite{he2023group}&Differential form&$\ell^{1,0}$-constraint&RR\\\hline
CaSLR~\cite{he2024much}&Differential form&$\ell^{1,0}$-constraint& RRC \\\hline
WG-IDENT~\cite{tang2025wgidentweakgroupidentification}&Weak form&$\ell^{1,0}$-constraint& RR\\\hline
\end{tabular}
\caption{A summary of identification methods of differential equation reviewed in this article.}\label{tab_summary}
\end{table}

\section{Coefficient recovery and  further  refinements}\label{sec:refinement}

In this section, we present some methods that help to further stabilize the identification process and  improve the identification accuracy. 

\subsection{High dynamic region and narrow fit}\label{sss:narrow}

During the feature selection process or after the support of the optimal model is determined, the coefficient values are computed via least-squares fitting of all the given data.  
However, the regions over which the differential equation  has larger variations seem to provide more accurate information characterizing the underlying differential equation, and this is explored in \cite{tang2023weakident}  by considering the  high dynamic region.  

Let $\{\Omega_h\}_{h=1}^H$ be a cover of the computational domain. In \cite{tang2023weakident}, $\Omega_h$ is the support of the $h$-th test function, and the feature term $u u_x$ is used to identify the high dynamic region, and $s(h)$ is defined to measure the variations of dynamics on each of $\Omega_h$.  In fact, this $s(h)$ can be computed by the leading coefficient error of $u u_x$, see \eqref{e: leading coefficients} in Section \ref{ssec:normalization}, here we set $k$ as the index of $u u_x$. 

Consider the set $\mathcal{S}=\{ s(h) | h=1, \dots, H\}$.  This set is divided into  mildly and highly dynamic regions automatically: partitioning the histogram of $\mathcal{S}$ into $N_{\cal S}$ bins $(b_1,b_2,...,b_{N_{\cal S}})$, we consider the cumulative sum of the bins $B(j) = \sum_{i=1}^j b_i$, then fit the  function $B(j)$ with a piecewise linear function $r(j)$ with one junction point $\Gamma$ using the cost function $ \sum_j (B(j)-r(j))^2/B(j)^2$.  The junction point $\Gamma$ separates the highly dynamic   and  mildly dynamic regions. 
Any $h$ with  $s(h) \geq \Gamma$ gives the high dynamic region  $\Omega_h$ that is used for the coefficient recovery.  Let the collection of the row indices of highly dynamic regions  be an ordered set:
\begin{equation*}
    \mathbb{H} =\{ h_i \;| \; \bar{s}(h_i) \geq \Gamma, \; h_i< h_j \text{ for } i <j \}.
\end{equation*}

We consider a submatrix  using only the ordered rows from  the highly dynamic region  $\mathbb{H}$ for both $F$ and $\bb$, i.e., 
$ F_{\rm narrow} := F^{\mathbb{H}}$ and $\bb_{\rm narrow}: = \bb^{\mathbb{H}}$. 
The coefficients are computed by solving 
$$
F_{\rm narrow}{\bc}=\bb_{\rm narrow}.
$$
This matrix is thus narrower than using the full data, which makes the computation easier, and gives good coefficient value recovery also.

\subsection{Feature matrix column normalization} \label{ssec:normalization}

The least-squares is used in many methods for coefficient value recovery for differential equation identification.  Since the accuracy of least squares is highly dependent on the conditioning of the feature matrix  \cite{bjorck1990least,bjorck1991error}, two types of normalization for the columns of the feature matrix is used to improve the coefficient recovery. 

The first type of normalization is \textbf{column normalization}, which is applied to the feature matrix as an input to SP, and  each column of $F$ is normalized by its own norm:
\begin{equation*} %\label{eq:Wdagger}
  F' = \displaystyle{ \left[\frac{F_1}{\|F_1\|}, \frac{F_2}{\|F_2\|}, \dots, \frac{F_{N_f}}{\|F_{N_f}\|}\right]}.
\end{equation*}
We observe that the scale of the columns in the feature matrix usually varies  substantially  from column to column, which negatively affects the SP.
This column normalization helps to prevent a large difference in the scale among the columns. 

The second normalization is \textbf{error normalization}, which is particularly effective for coefficient value recovery. 
For Weak-IDENT, the columns in $F$ are given by certain  derivatives of a monomial of $u$. When we compute the feature matrix with noisy data, the noise has different effects on different features. Dentoe the $k$-th feature as $ \frac{\partial^{\alpha}}{\partial x ^{\alpha}} \left( u^{\beta}\right)$. Over each $\Omega_h$,
the noisy data with noise $\epsilon$ in   \eqref{eq:givedata} give rise to the leading coefficient error:
\begin{equation}
s(h,k) = \beta \left| \int_{\Omega_h} u^{\beta -1} \frac{\partial^{\alpha}}{\partial x ^{\alpha}}   \left(\phi_h(x,t)\right) dxdt \right|, \;\; h = 1,2,..., H, \; \beta \ge 1.
\label{e: leading coefficients}
\end{equation}
When $\alpha = \beta = 0$, we set $s(h,k)=1$. This leading coefficient $s(h,k)$ depends on the row index $h$ and the column index $k$. For the $k^{\rm th}$ column, we define
$$\langle s(h,k)\rangle_{h} = \frac{1}{H} \sum_{h=1}^H s(h,k)$$
as an average of these leading coefficients over the rows, and normalize $F$ as \begin{equation} \label{eq:Wtilde}
  \tilde{F} = \displaystyle{ \left[\frac{F_1}{\langle s(h,1)\rangle_{h}}, \frac{F_2}{\langle s(h,2)\rangle_{h}}, \dots, \frac{F_{N_f}}{\langle s(h,{N_f})\rangle_{h} }\right].}
\end{equation}
For the $h$-th test function, we let $s(h)$ be the average leading coefficient error over selected features, especially $uu_x$ in Weak-IDENT.
This error normalization not only helps with more uniform entry values with less variation across different columns, but also used in finding the high dynamic regions as in the previous section.

\section{Results of identifying constant coefficient differential equation}
\label{sec:numerical}

We summarize identification results of the constant coefficient differential equations using the methods described above. 

In IDENT \cite{kang2021ident}, results are presented for Burgers' equation (15) with 0 to 30 \% noise, 
\[
u_t + u u_x= 0, \quad x \in [0,1], \quad
u(x,0)= \sin 4\pi x, \quad u(0,t)= u(1,t)= 0, 
\]
and viscous Burgers' equation (16), $ u_t + u u_x= 0.1 u_{xx}$ with the same initial and boundary condition with 0 to 0.12 \% noise. 
We say that the noise is $p\%$ by setting 
\begin{equation}\label{eq:sigma_ident}
\sigma=\frac{p}{100}\sqrt{\frac{1}{N_xN_t}\sum_n \sum_{\mathbf{i}} \left(u(x_i,t_n)\right)^2} .
\end{equation}
A new noise-to-signal ratio is proposed: 
\begin{equation}
\label{eqnsr}
\text{Noise-to-Signal Ratio (NSR)} := \frac{\|{F} \mathbf{c} -\bb\|_{L^2}}{\min_{j:\ c_j\neq 0} \|{F}[j]\|_{L^2} | c_j| }. 
\end{equation}
This is derived from \cite[Theorem 1]{kang2021ident}, showing that the signal level is contributed by the minimum of the product of the coefficient and the column norm in the feature matrix, i.e., $\min_{j:\ c_j\neq 0} \|{F}[j]\|_{L^2} | c_j|$.  This is used to explain the difficulties of identification of some differential equations.   Table \ref{table:nsr} shows NSR for Burgers' equation and viscous Burgers' equation.  While NSR values are similar, this  corresponds to big difference in noise level represented as \% in \eqref{eq:givedata}.  In the case of viscous burger's equation, even if it seems IDENT can only recover up to 0.12\% of noise, this corresponds to much more difficult problem compared to Burgers' equation, since NSR is 12.13 vs.~3.23.  IDENT shows good results for low level of noise, and a low order of underlying differential equations.  Comparisons are presented against Robust IDENT \cite{he2022robust} and Weak IDENT \cite{tang2023weakident} in these papers. 
\begin{table}
\centering
$u_t = - u u_x$  with $u(x_i, t_n) + \epsilon_i^n$

\begin{tabular}{c|cccccccc}
\hline
$\epsilon$ noise \% & 0 & 4 & 8 & 12 & 16 & 20 & 24 & 28\\
NSR &
0.04 & \textbf{2.18} & 3.09 & 3.33 & 3.40 & 3.31 & 3.31 & 3.23\\
\hline
\end{tabular}

\vspace{0.3cm}
$u_t = - u u_x + 0.1 u_{xx}$  with $u(x_i, t_n) + \epsilon_i^n$ 
\begin{tabular}{c|cccccccc}
\hline 
$\epsilon$ noise  \% & 0 & 0.02 & 0.04 & 0.06 & 0.08 & 0.10 & 0.12 \\
NSR & 0.05 & \textbf{2.02} & 4.04 & 6.06 & 8.08 & 10.10 & 12.13 \\
\hline 
\end{tabular}
\vspace{0.3cm}
\caption{[IDENT] Difficulties of LASSO identification represented in terms of new Noise-to-Signal Ratio (NSR) \eqref{eqnsr}. For similar values of NSR (in bold), while this corresponds to 4\% noise in Burgers' equation, it corresponds to 0.02 \% noise for viscous Burgers' equation.  These are equations (15) and (16) in \cite{kang2021ident}. This represents the difficulty in identification for viscous Burgers' equation in this case for a small level of \% noise. }  \label{table:nsr}
\end{table}

In Robust IDENT \cite{he2022robust},  various measures are introduce to systematically compare the identification results. Let $\bc^*$ be the true coefficient vector and $\bc$ be the recovered coefficient vector.
The relative coefficient error $e_c$ 
\begin{equation}\label{E:coeffError}
e_c=\frac{\|{\mathbf{c}-\mathbf{c}^*\|_1}}{\|\mathbf{c}^*\|_1},
\end{equation}
measures the accuracy in the recovery of PDE coefficients, the grid-dependent residual error $e_r$ 
\begin{equation}
e_r=
\begin{cases}
\sqrt{\Delta x \Delta t}\|F({\mathbf{c}}-\mathbf{c}^*)\|_2 \mbox{ for 1D PDE},\\
\sqrt{\Delta x \Delta y\Delta t}\|F({\mathbf{c}}-\mathbf{c}^*)\|_2 \mbox{ for 2D PDE},
\end{cases}\; 
\label{E:resError}
\end{equation}	
measures the difference between the learned dynamics and the denoised one by SDD, and the evolution error (dynamic error)
\begin{align}\label{E:dynError}
  e_e=\Delta x \Delta t\left(\sum_n\sum_{\mathbf{i}} |u({x}_{i},t_n)-\hat{u}({x}_{i},t_n)|\right)
\end{align}
where $u$ and $\hat{u}$ denote the solution of the exact and identified PDE from the same initial condition, respectively, measures  how well the solution of the identified PDE matches the dynamics of the correct PDE.
In Weak IDENT \cite{tang2023weakident}, for the coefficient error \eqref{E:coeffError}, $\ell^2$ norm $E_2$ as well as $\ell^\infty$ norm $E_{\infty}$ are considered, with residual error \eqref{E:resError} as $E_{res}$ and dynamic error \eqref{E:dynError} as $E_{dyn}$ in addition to TPR and PPV.  The True Positive Rate (TPR) measures the fraction of features that are found out of all features in the true equation, and is defined as the ratio of the cardinality of the correctly identified support over the cardinality of the true support:
\begin{equation*}
\text{TPR} = | \{l: {c}^*_l \neq 0, \; {c}_l \neq 0 \}| / | \{ l: {c}^*_l \neq 0\}|.
%\label{E:tprError}
\end{equation*}
The TPR is 1 if all the true features are found. The Positive Predictive Value (PPV) indicates the presence of false positives: 
\begin{equation*}
\text{PPV} = | \{l: {c}^*_l \neq 0, \; {c}_l \neq 0 \}| / | \{ l: {c}_l \neq 0\}|, 
%\label{E:ppvError}
\end{equation*}
it is the ratio of the cardinality of the correctly identified support over the total cardinality of the identified support.

\begin{table}
\centering
\begin{tabular}{c|c|cc}
\hline
& Identified equations & $e_c$ & $e_r$\\
\hline 
\cite{kang2021ident} &  $u_t = - 1.02  u u_x$  & $2.39 \times 10^{-2}$ & $8.27 \times 10^{-5}$ \\
\cite{he2022robust} &  $u_t = - 1.02  u u_x$  & $2.39 \times 10^{-2}$ & $8.27 \times 10^{-5}$ \\
\cite{rudy2017data} &  $u_t = - 0.38  u u_x - 0.93 u$  & $1.56$ & $1.84 \times 10^{-3}$ \\
\cite{schaeffer2017learning} & $u_t = -0.02 u u_x + 1.84 u^2 + 10.67 u -1$ & $13.59$ & $7.16 \times 10^{-3}$ \\ \hline 
\end{tabular}
\caption{[Robust IDENT] Comparison of identification of Burgers' equation $u_t = - u u_x$ with initial  $u(x, 0) = \sin( 4\pi x) \cos (2 \pi x)$, the equation (4.7) in  \cite{he2022robust} for 40 \% of  noise as in \eqref{eq:givedata}.  
IDENT \cite{kang2021ident} and Robust IDENT \cite{he2022robust} gives identical results in the first two rows. }  \label{table:robustCompar}
\end{table}
In Robust IDENT \cite{he2022robust},
Burgers' equation $u_t = - u u_x$ is considered for 0 to 40 \% noise and it identifies the result to be $u_t = -0.7366 u u_x$ for 40 \% noise with $e_c = 2.63 \times 10^{-1}$ and $e_r = 1.64 \times 10^{-1}$.  
This is compared to \cite{rudy2017data} and \cite{schaeffer2017learning} for noise level 0 to 40 \% in Table 5 in \cite{he2022robust}. For a fair comparison initial condition $u(x, 0) = \sin( 4\pi x) \cos (2 \pi x)$ is used, and for 40 \% of noise, IDENT and Robust IDENT show good results which is  summarized in Table \ref{table:robustCompar}.
Viscous Burgers' equation $u_t= - u u_x + 0.1 u_{xx}$ is considered for 0 to 5 \% noise.  Compared to \cite{kang2021ident}, Robust IDENT shows more stability handling up to 5\% noise vs. 0.12 \%.  For 5 \% noise, it identifies $ u_t = -1.0170 u u_x + 0.0976 u_{xx}$ with $e_c = 1.77 \times 10^{-2}$ and $e_r = 2.21 \times 10^{-2}$.  In addition, Transport equation $u_t = -u_x$ is considered for both smooth and discontinuous initial conditions.   It  identifies the result to be $u_t = -0.9421 u_x$ for 30 \% noise with $e_c = 5.79 \times 10^{-2}$ and $e_r = 4.31 \times 10^{-2}$.  
KdV equation $ u_t+ 6 u u_x + u_{xxx} =0 $ is explored, identifying the equation to be $ u_t+ 6.135 u u_x + 1.0580  u_{xxx} =0$ with $e_c = 2.77 \times 10^{-2}$ and $e_r = 1.21$.   In addition, larger dictionary, 2-dimensional PDE and identifiability based on the given data is explored in \cite{he2022robust}.

\begin{table}[t]
\centering
\begin{tabular}{l|l|c}
\hline
Method & Equation considered & Tested noise level  \\ \hline
IDENT & $u_t = - u u_x$ & $p=0$ to $30 \%$ \\
 & $u_t = - u u_x + 0.1 u_{xx}$ &  $p=0$ to $0.12 \%$ \\ \hline
Robust IDENT & $u_t = - u u_x$ & $p= 0$ to $40 \% $\\
& ($u_t = -0. 7366 u u_x$) & $(p = 40 \%)$ \\ 
 & $u_t = - u u_x + 0.1 u_{xx}$ & $p=0$ to $5 \%$  \\ 
 & ($u_t = - 1.0170 u u_x + 0.0976 u_{xx}$) & $(p = 5 \%)$ \\
& $u_t = -u_x$ & $p= 0$ to $30 \% $ \\
& ($u_t = - 0.9421 u_x$) & ($p=30 \%$) \\
& KdV: $ u_t+ 6 u u_x + u_{xxx} =0 $ & no noise  \\
& ($ u_t+ 6.135 u u_x + 1.0580  u_{xxx} =0$) & \\ \hline
Weak IDENT & $u_t = - u_x + 0.05 u_{xx}$ & $\sigma_{\rm NSR}=$ 0 to 1 \\
& ($u_t = - 1.00792 u_x + 0.05029 u_{xx}$) & ($\sigma_{\rm NSR}=$ 1) \\
& $u_t  =  0.30 (u^2)_{yy}
- 0.80 (u^2)_{xy} + 1.00 (u^2)_{xx}$ &$\sigma_{\rm NSR}=$ 0 to 0.1\\
& ($ u_t  =  0.29 (u^2)_{yy} -0.79 (u^2)_{xy} +0.99 (u^2)_{xx}$) 
     & ($\sigma_{\rm NSR}=$ 0.08) \\
& $u_t  =  + v^3 + u +0.1 u_{yy} +0.1u_{xx} - uv^2 +u^2v - u^3$ & $\sigma_{\rm NSR}=$ 0.01 to 0.1\\
& $v_t  =   v +0.1v_{yy} +0.1v_{xx} - v^3 -uv^2 - u^2v -u^3$ & \\
& KdV:  $ u_t + 0.5 u u_x + u_{xxx} =0 $& $\sigma_{\rm NSR}=$ 0 to 0.2 \\
& KS: $ u_t + u + u_{xx} + u_{xxxx} =0 $ & $\sigma_{\rm NSR}=$ 0 to 1 \\
%& NLS & $\sigma_{\rm NSR}=$ 0 to 0.5 \\ 
\hline
\end{tabular}
\caption{Some experiments presented in IDENT \cite{kang2021ident}, Robust IDENT \cite{he2022robust} and Weak IDENT \cite{tang2023weakident}. The rows with parenthesis shows one representative result for each method. Compared to IDENT, Robust IDENT can handle higher level of noise as well as higher derivative as $u_{xxx}$.  Comparied to Robust IDENT, Weak IDENT can handle even higher order differential equations with higher level of noise.  With the new definition, $\sigma_{\rm NSR}=1$ represents 100 \% of noise.}
\label{table:summary}
\end{table}

In Weak-IDENT \cite{tang2023weakident}, more extensive experiments and comparisons are presented. Gaussian noise is added with $\epsilon_i^n \sim \mathcal{N}(0, \sigma)$ for $\epsilon_i^n$, and $U_i^n$ in \eqref{eq:givedata}.
We note that the noise level is redefined considering the values of $U$:  
\begin{equation} 
%\epsilon_i^n\sim \mathcal{N}(0,\sigma)~\text{where}~\sigma = 
\sigma = \sigma_{\text{NSR}}\sqrt{\frac{1}{N_tN_x}\displaystyle \sum_{{i},n}|u(x_i, t_n) - (\max_{i, n} u(x_i, t_n) +\min_{i, n} u(x_i, t_n) )/2 |^2}
\label{E:newNoise}
\end{equation} 
for $i = 1,...,N_x$, $n = 1,...,N_t$.  This definition of noise level reflects the local variation of the given data.  When $p/100=\sigma_{\text{NSR}}$, the standard deviation of the noise in~\eqref{E:newNoise} tends to be smaller than that in~\eqref{eq:sigma_ident}.  
Weak IDENT \cite{tang2023weakident} presents PDE experiments and comparisons with IDENT \cite{kang2021ident}, Robust IDENT \cite{he2022robust}, and methods using the weak form such as \cite{reinbold2020using} and Weak SINDy for PDE \cite{messenger2021weak} for transport equation with diffusion, Korteweg-de Vires (KdV), Kuramoto-Sivashinsky (KS), 1-dimensional system of Nonlinear Schr\"{o}dinger equation (NLS), 2-dimensional Anisotropic Porous Medium equation (PM) and 2-dimensional system of Reaction-Diffusion equation.  Statistic experiments are presented with 50 experiments for each level of $\sigma_{\rm NSR}\in \{ 0.01,0.1, 0.2, ,...,0.9\}$.  Summary of some equations experimented in IDENT \cite{kang2021ident}, Robust IDENT \cite{he2022robust} and Weak IDENT \cite{tang2023weakident} is presented in Table \ref{table:summary}.   
For periodic boundary condition, in \cite{tang2023fourier} boundary padding method is proposed to handle non-periodic boundary conditions. 

In \cite{tang2023weakident}, additional comparison is presented for dynamical systems and compared against SINDy \cite{brunton2016discovering} and Weak SINDy for dynamical systems \cite{messenger2021weakode} for $\sigma_{\rm NSR} =0.1$ and $0.2$. Dynamical system considered are 2-dimensional linear system, 2-dimensional nonlinear Van der Pol equation, 2-dimensional nonlinear Duffing equation, 2-dimensional nonlinear Lotka-Volterra equation, and 3-dimensional nonlinear Lorenz equations.  In \cite{kaptanoglu2023benchmarking}, various methods are compared and Weak SINDy showed good performances.  It is shown in \cite{tang2023weakident} that Weak IDENT outperforms or comparable to Weak SINDy \cite{messenger2021weak, messenger2021weakode}, and that SP with narrow-fit given by high dynamic region and trimming stabilizes the identification process.

\section{Identification of differential equation with varying coefficients}\label{sec:varying}

The methods discussed in the previous sections consider identifying PDEs with constant coefficients.  In many real physical phenomena, one cannot assume that the dynamics is given by one equation, but the dynamics may change with time.    We present methods to identify differential equations with varying coefficients. 

%\subsection{Identification of differential equation with varying coefficients}\label{ssec:varying}

\subsection{IDENT: group LASSO with BEE}\label{sec:GPLASSO}

In IDENT \cite{kang2021ident}, the varying coefficients are  modeled by finite element bases.  Let $\{ \phi_m\}_{m=1}^{N_b} $ be a set of $N_b$  number of finite element bases in the spatial domain, then the coefficients $c_k$ in \eqref{E:general_constant} can be represented as
\begin{equation}\label{E:L}
c_k (x) \approx \sum_{m=1}^{N_b} c_{k,m} \phi_m(x) \text{ for } k=1,\dots,N_f.
\end{equation}
Each feature $f_k$ corresponds to a set of coefficients $\{c_{k,m}\}_{m=1}^{N_b}$. Suppose that the time domain is discretized into $N_t$ grids. Denote the value of $f_k$ at $(x_i,t_n)$ by $(f_k)_i^n$, and $\phi_m(x_i)$ by $(\phi_m)_i$. The feature matrix $F$ is of size $(N_x N_t) \times (N_b N_f)$,
{\footnotesize
\begin{align} \label{E:general_F}
F = \left[ \begin{array}{ccc|c|ccc}
\vdots &  & \vdots &  &\vdots &  & \vdots  \\
(f_1)_i^n(\phi_1)_i & \dots & (f_{1})_i^n(\phi_{N_b})_i & \cdots &(f_{N_f})_i^n(\phi_1)_i & \dots  & (f_{N_f})_i^n(\phi_{N_b})_i \\
\vdots &  & \vdots & &  \vdots &  & \vdots  \\
\end{array}  \right],
\end{align}
}
the coefficient vector to be identified is 
\[ 
\mathbf{c} = \left[c_{1,1},\cdots,c_{1,N_b} | c_{2,1}, \cdots, c_{2,N_b} | \cdots  | c_{N_f,1}, \cdots, c_{N_f,N_b} \right]^{\top}
\in \mathbb{R}^{N_b N_f},  
\]
and the feature response is
$$
\mathbf{b}=\left[u_t(x_1,t_1),u_t(x_2,t_1),\cdots,u_t(x_{N_x-1},t_{N_t}), u_t(x_{N_x},t_{N_t})  \right]^{\top} \in \mathbb{R}^{N_xN_t}.
$$
The feature matrix $F$ has a block structure that its columns can be grouped into $N_f$ groups, each of which corresponds to a feature, with $N_b$ columns in each group.

To identify the varying coefficients in \cite{kang2021ident}, group LASSO \cite{yuan2006model} is used, and a set of possible active features is selected by 
solving 
\begin{equation}\label{eq:groupident}
    \widehat{\mathbf{c}}(\lambda) ={\textstyle \arg \min_{\bc}} \left\{ \frac{1}{2} \| \bb -  F \mathbf{c} \|^2_2 + \lambda \sum_{k=1}^{N_f} \left(\sum_{m=1}^{N_b} |c_{k,m}|^2 \right)^{\frac{1}{2}} \right\},
\end{equation}
where ${F}$ is the feature matrix in \eqref{E:general_F} and ${\mathbf{b}} $ is the feature response computed from given data.
In (\ref{eq:groupident}), the second term on the right-hand side is a penalty on the group sparsity of $\bc$ of $\ell^{2,1}$ type regularization: it computes the $\ell^2$ norm of each group of $\bc$, then puts them together, and computes the $\ell^1$ norm.

In \cite{kang2021ident}, Base Element Expansion (BEE) is proposed to select the number of bases $N_b$. 
For each fixed $N_b$, group LASSO gives candidate coefficient vector $\widehat{\bc}(\lambda)$, and the normalized block magnitudes from group LASSO is recorded as  
\[
\mathcal{B}(N_b) :=  
 \left\{\|{F}[k]\|_{1}\left\| \sum_{m=1}^{N_b} \frac{\widehat{\bc}(\lambda)_{k,m}}{\|{F}[k,m]\|_\infty} \phi_m \right\|_{L^1}\right\}_{k=1,\ldots,N_f},
\]
where ${F}[k,m]$  denotes the $m$-th column (element) of the $k$-th group of ${F}$, and ${F}[k]$ denotes the $k$-th group.
Following the idea of the convergence of the finite element approximation, as more basis functions are used, i.e. as $N_b$ increases, the more accurate the approximation will be.  In the BEE procedure, the normalized block magnitudes reach a plateau as $N_b$ increases, and the number of bases is chosen as the smallest $N_b$ that $\mathcal{B}(N_b)$ reaches a plateau.  Some experiments show stabilization with $N_b=5$, and $N_b=20$ was typically chosen.  With a properly chosen $N_b$, TEE is applied to find the best identification in IDENT \cite{kang2021ident}. 

\subsection{GP-IDENT: $\ell^{1,0}$-group sparsity constraint}
\label{sec:GPIdent}

% Using subspace pursuit shows more stability compared to using $\ell^1$ LASSO. 
Group SP is proposed in \cite{he2023group} for $\ell^{1,0}$-group sparsity constraint.   
In \cite{he2023group}, the construction of linear system is similar to \eqref{E:L}, but using B-spline. The B-spline bases are defined over space-time domain and each basis is a product of a B-spline in space and another B-spline in time.
% : 
% \[
% c_k(x_i,t_n)f_k(x_i,t_n)\approx \sum_{m=1}^Mc_{k,m} B_m(x_i,t_n)f_k(x_i,t_n). 
% \]
GP-IDENT proposes to generate a set of candidate PDEs by considering 
\begin{align}
\widehat{\mathbf{c}}=\argmin_{\mathbf{c}}&\|{F}\mathbf{c}-{\mathbf{b}}\|^2_2~\mbox{ subject to }\|\mathbf{c}\|_{\ell^{1,0}} = k, \label{eq:gpident}
\end{align}
where the $\ell^{1,0}$-norm of a vector $\bc\in\mathbb{R}^{N_b N_f}$ is
\begin{align*}
\|\mathbf{c}\|_{\ell^{1,0}} :=\left\|\begin{bmatrix}\|\mathbf{c}_1\|_1 & \dots &\|\mathbf{c}_{N_f}\|_1\end{bmatrix}\right\|_0
\end{align*}
which counts the number of groups with non-zero coefficients. Here the notation $\mathbf{c}_k = \left[c_{k,1},\cdots,c_{k,N_b}  \right]^{\top}$ is used.
The constraint enforces group sparsity by explicitly specifying that only $k$ groups of features have nonzero coefficients. The solution of~\eqref{eq:gpident} corresponds to a PDE model with exactly $k$ features that best fits the given data.

In \cite{he2023group}, Group Projected Subspace Pursuit (GPSP) is proposed.  It is a greedy algorithm that refines the selected features by iterations. Denote the $k$-th group of $F$ by $F_g$. For a fixed level of group sparsity $k^*\geq 1$,  suppose the set of group indices selected by the $(l-1)$-th iteration is $T^{l-1}$, and define the residual of fitting the data  using groups specified by indices in $T^{l-1}$ as
\begin{align*}
\mathbf{\by}^{l-1}_r =\resid(\by, F_{T^{l-1}})= \mathbf{y}-\text{proj}(\by,F_{T^{l-1}})=\by-F_{T^{l-1}}F_{T^{l-1}}^\dagger\mathbf{y} .%\label{eq_y_residual}
\end{align*} 
 Here $F_{T^{l-1}}$ is obtained by concatenating the group features $\{F_k\}_{k\in T^{l-1}}$ horizontally. The proposed scheme consists of two stages in each iteration. (1) Expand $T^{l-1}$ to $\wT^{l}$: For the $l$-th iteration, we take the union of $T^{l-1}$ with the set of $k^*$ groups with the highest $k^*$ values of 
\begin{align*}
P(\by_r^{l-1},F_g)=\frac{\left|\text{proj}(\by_r^{l-1},F_k)^{\top}\by_r^{l-1}\right|}{\|\text{proj}(\by_r^{l-1},F_k)\|_2\|\by_r^{l-1}\|_2} %\label{eq:PyrPg}
 \end{align*}
and denote the union set as $\wT^{l}$. 
(2) Shrink $\wT^{l}$ to $T^l$.
 Let $\mathbf{x}^l_p = F_{\wT^l}^\dagger\by$.
We project $\by$ to the column space of $F_{\wT^l}$ with decomposition $$\by_p = \text{proj}(\by,F_{\wT^l})= \sum_{k\in \wT^l }F_k \mathbf{x}^l_p[k],$$ where $\bx^l_p[k]$ is the subvector of $\bx^l_p$ corresponding to the $k$-th group.
For $k\in\wT^l$, its norm $\|F_k \mathbf{x}^l_p[k]\|_2$ provides a measure of the importance of the $k$-th group, that we keep indices of $k$ most important groups and remove the others. The refined set of indices is $T^{l}$. 
After the $l$-th iteration, we compute $\by_{r}^{l} =\resid(\by, F_{T^{l}})$. If $\|\by_{r}^{l}\|_2 > \|\by_r^{l-1}\|_2$, we set $T^{l} = T^{l-1}$ and take the groups indexed by $T^{l}$ as our final selection;  otherwise, repeat the procedure.

The above described Group Projected Subspace Pursuit (GPSP) was proved to converge under suitable conditions in  \cite{he2025group}. 
Specifically, if the Block Restricted Isometry Property (BRIP) condition is satisfied, \cite[Theorem 1]{he2025group} proves the convergence of GPSP for clean data.  When the measurements contain an additive perturbation, \cite[Theorem 2]{he2025group} gives the characterization of the recovery distortion in terms of the BRIC of the sampling matrix and the magnitude of the perturbation. 
As demonstrated in \cite{he2023group,he2025group}, compared to other greedy algorithms for group sparsity, GPSP has more robust performances in general applications.

%\subsubsection{Model selection for group sparsity}

For the model selection, RR \eqref{eq:RR} is used in \cite{he2023group}, and in \cite{he2024much}, the Reduction in Residual Curve (RRC) \eqref{eq_modelscore} is proposed. These methods are particularly useful for varying coefficient PDEs, as their feature matrices have larger sizes.  

In addition, Weak formulation based GP-IDENT is investigated in WG-IDENT \cite{tang2025wgidentweakgroupidentification} combining the benefits of both Weak IDENT with GPSP. 

\subsubsection{CaSLR: patch-based local and global regession}\label{sec:CaSLR}
In  large-scale modeling such as  meteorology, data are typically measured by local sensors distributed across scattered locations. Assuming the existence of a global governing equation with spatially varying coefficients, the authors of~\cite{he2024much} proposed the Consistent and Sparse Local Regression (CaSLR) method. This patch-based identification approach combines local regression with global consistency constraints enforced through group sparsity.

% Assume that the unknown PDE takes the following form:
% \begin{align}
% u_t(x,t) = \sum_{k=1}^Kc_{k}(x,t)f_k(x,t),\label{eq_PDEform}
% \end{align}
% where $\mathcal{F} = \{f_k:\Omega\times [0,T]\to\mathbb{R}\}_{k=1}^K$ is a dictionary of features that contains partial derivatives of $u$ with respect to the space (e.g., $u_x,u_{xxx}$, etc.), functions of these terms (e.g., $uu_x$ and $u^2$, $\sin (u)$ etc.); and  $c_k:\Omega\times [0,T]\to\mathbb{R}$, $k=1,2,\dots, K$ are the respective coefficients.  We assume that the dictionary is rich enough that it is over-complete, i.e., the underlying PDE is expressed as~\eqref{eq_PDEform} when at least one of the feature's coefficients is null. One also has to assume the measurements can resolve local variation of the coefficients, i.e., 
Let $\{f_k\}_{k=1}^{N_f}$ be a set of candidate features, and let $\{\Omega_j\}_{j=1}^J$ be subsets of the spatiotemporal observation domain. If the coefficients of the underlying PDE vary smoothly across the domain, the dynamics within each patch $\Omega_j$ can be well approximated by a PDE with constant coefficients. For each patch $\Omega_j$, let $\widehat{\mathbf{c}}_j = (\widehat{c}_1^j, \ldots, \widehat{c}_K^j)$ denote the estimated coefficient vector derived from local data. The local regression error for $\Omega_j$ is given by
\begin{align*}
\mathcal{E}_{\text{loc}}^j(\widehat{\mathbf{c}}_j) &= \int_{\Omega_j} \left( u_t(x,t) - \sum_{k=1}^{N_f} \widehat{c}_k^j f_k(x,t) \right)^2 \,dx\,dt,  
%\label{eq_loc}
\end{align*}
and the total global regression error aggregates these local errors:
\begin{equation}
\mathcal{E}(\widehat{\mathbf{c}}) = \sum_{j=1}^J \mathcal{E}_{\text{loc}}^j(\widehat{\mathbf{c}}_j), \quad \text{where} \quad \widehat{\mathbf{c}} = [\widehat{c}^1_1, \widehat{c}^2_1,\cdots, \widehat{c}^J_1,\widehat{c}^1_2,\cdots,\widehat{c}^J_{N_f}].
\label{eq_global}
\end{equation}
Here $\widehat{\bc}$ consists of $N_f$ groups, and each group includes the coefficients of the same feature across different domains.
To ensure that the coefficients identified from all patches have the same support,  the $\ell^{1,0}$-constrained group sparsity is used:
\begin{equation*}%\label{eqCaSLR}
\widehat{\mathbf{c}}^{l} = \arg\min_{\widehat{\mathbf{c}}}\mathcal{E}(\widehat{\mathbf{c}}) \quad
\text{subject to: }\quad \|\widehat{\mathbf{c}}\|_{\ell^{1,0}} = l.
\end{equation*}
Here the group structure is organized so that
\begin{equation}\nonumber
\|\widehat{\mathbf{c}}\|_{\ell^0} = \|[\|\widetilde{\mathbf{c}}_1\|_1,\dots,\|\widetilde{\mathbf{c}}_{N_f}\|_1]\|_0,\\
\end{equation}
where $\widetilde{\mathbf{c}}_k=[\widehat{c}_k^1,\dots,\widehat{c}_k^J]\in\mathbb{R}^J, k=1,2,\dots,N_f$.
Setting $\widehat{\mathbf{c}}^0 = [0,\dots,0]$, the Residual Reduction Curve (RRC) is proposed, and the score for the $l$-th model is computed by
\begin{align}
S^{l} =  \mathcal{E}(\widehat{\mathbf{c}}^{l})+\rho\frac{l}{N_f}\label{eq_modelscore}
\end{align}
for $l=1,2,\dots, N_f-1$, where $\rho>0$ is a penalty parameter for using more terms from the dictionary and $\mathcal{E}$ given in \eqref{eq_global}.
The candidate with $l^*$ features is defined to be the optimal if $S^{l^*}=\min_{l=1,\dots,N_f-1}S^l$. 
The metric $S^k$ in~\eqref{eq_modelscore} evaluates a candidate model with $l$ features by considering two factors: the fitting error and the model complexity penalty controlled by the parameter $\rho$.  In \cite{he2025group},  $\rho$ is fixed to be the mean of 
$\{\mathcal{E} (\widehat{\mathbf{c}}^{k})\}_{k=1}^{N_f}$ to balance the two components  in \eqref{eq_modelscore}.

The above local approximation works exactly when: 1) the variable coefficients are bounded away from zero and vary slowly on the patch, 2) the solution data contain diverse information content, on the patch. The first condition basically requires the mathematical problem to be well-posed and the second condition says the local regression system determined by the solution data is well-conditioned. 

\subsection{Numerical results of identification of differential equation with varying coefficients}

In IDENT \cite{kang2021ident}, there are three sets of experiments. 
In the first set, it is assumed that one coefficient is known a priori to vary in spacial dimension. The experiment was on viscous Burgers' equation where the diffusion term is assumed to have a varying coefficient.  It was experimented with no noise and with 2 \% noise as in \eqref{eq:givedata}, and with or without LSMA denoising.  
For the second set of experiments, two coefficients are assumed to be known a priori to vary in spacial direction for advection-diffusion equation.  In the final experiment, all coefficients are free to vary without any a priori information.  IDENT with group LASSO, BEE and TEE showed good results based on the ideas of numerical convergence. 

In GP-IDENT \cite{he2023group}, more extensive experiments with comparison is provided.  In addition to coefficient error, Jaccard index \cite{jaccard1912distribution} is used to quantify the coefficient support identification accuracy
\begin{align*}
J(\widehat{T},T^*) =\frac{|\widehat{T}\cap T^*|}{|\widehat{T}\cup T^*|},
\end{align*}
where $\widehat{T}$ denotes the group index set  in the identified model,  $T^*$ is the group index set in the true equation, and $|\cdot|$ gives the number of elements in the set.  Note that $J(\widehat{T},T^*)=1$ if and only if $\widehat{T}=T^*$, i.e., the underlying model is exactly identified. 

The differential equation considered in GP-IDENT \cite{he2023group} are Advection-diffusion equation, Fisher's equation, viscous Burgers' equation, KdV equation, KS equation, Schr\"{o}dinger and nonliner Schr\"{o}dinger equation for space and time varying coefficients and 0 to 2\% noise as in \eqref{eq:givedata}.
The Group Projected Subspace Pursuit (GPSP) is compared with Block Subspace Pursuit \cite{kamali2013block} in detail.  One of the experiments shows the stability against the size of dictionary where the size of dictionary is ranging from $N_f = 35$ to $330$ feature terms.  
Compared to \cite{yuan2006model}, parametric partial differential equations identification \cite{Rudy2019DataDrivenIO}, and low rank recovery of PDE \cite{li2020robust}, Block SP and GPSP shows the best stable identification, correctly identifying feature terms, i.e., correct coefficient support even with a large dictionary.  GP-IDENT consistently produces accurate and robust results.

In addition, we note that WG-IDENT \cite{tang2025wgidentweakgroupidentification} explores using weak form with GP-IDENT, which gives stable recovery.

\section{Concluding remarks}\label{sec:summary}

This series of works propose methods for identifying differential equations using explainable, mathematically motivated techniques based on numerical PDE methods, convergence analysis, and least-squares fitting. Without requiring any training, we explore approaches to derive the best approximation of the underlying equation, even from  one single noisy trajectory of the solution. Since these works are fundamentally driven by the goal of understanding real-world phenomena, it is essential to employ interpretable methods. These techniques provide optimally fitted results for datasets of any size.

We explore and study each step carefully and provide mathematical frameworks for stable identification. 
Just by assuming that the underlying differential equation can be expressed as a linear combination of linear and nonlinear differential terms, this inverse problem of identification can be formulated as a simple linear system.  This formulation includes many well-known equations of much interest such as Korteweg–De Vries (KdV) equation, Kuramoto–Sivashinsky (KS) equation, reaction-diffusion equation, nonliner Schr\"{o}dinger equation and others. 

The foundation of identification is  model-selected least-square fitting, which contributes to finding the support as well as the values of the coefficient vector.  Firstly, it is important for the linear system to be well-conditioned.   Good denoising and derivative approximations for the feature matrix are important for stable identification.  We present denoising methods for finite difference methods as well as the effects of using the weak form.   SDD can handle more range of feature terms, while using test function and weak form shows stability for particular types of feature terms.  

Secondly, for least-square fitting, the more terms one has to fit, the smaller the residual will be, i.e. it will prefer to find more feature terms in identification.  Using sparsity is beneficial to prevent finding the most complex equations for this reason.  However, what feature terms to be identified should be considered carefully.  We propose mathematically motivated model selection methods, such as TEE, MTEE, BEE, etc., to explicitly identify the best equation.  As for sparse regression, we use SP, since SP gives $k$ feature terms with the approximately least residual in fitting $\| F\bc - \bb\|^2$ for each sparsity $k$, so that the result is  explainable. 

In addition, high dynamic region and narrow fit improve the least-square fit in coefficient value identification by using only the good data (high dynamic region), while ignoring many redundant unimportant data.  Trimming is also introduced to remove unnecessary feature terms that do not contribute to $u_t$ much.  This process is especially helpful since we use a large generic
dictionary.  % \roy{The following sentence is out of context of this paragraph; consider remove: We typically include all the feature terms defined by the number of multiples allowed as well as the order of derivative, e.g. for up to 2 multiples with second order derivative gives  a dictionary with  10 terms $\{1, u, u^2, u_x, u u_x, u_x^2, u_{xx}, u u_{xx}, u_x u_{xx}, u_{xx}^2 \}$, this is how we experimented with $N_f=330$. }  
Both methods contribute to having fewer data or feature terms with good information in least-square fitting.  

In summary, when only one single realization of time-dependent data is given, even if the observation is perturbed by noise, we propose methods to robustly  identify the underlying differential equation.  First, we set up a linear system with  denoising and good approximations to the derivatives in feature terms. Secondly, we 
consider reducing the number of feature terms by sparse regularization or constraint where we used SP for each sparsity level $k =1, \dots N_f$. This is then refined by mathematically motivated model selection methods, such as TEE, MTEE, CEE and others, which helps to pick an explainable result, e.g., best dynamic fit or best residual reduction.   
It will be exciting to apply these methods to real data which remains to be challenging due to the difficulties in the data acquisition. 

\section{Acknowledgement}

We would like to thank our co-authors of the papers reviewed in this paper: Prof.~Yingjie Liu for developing IDENT \cite{kang2021ident} together, Dr.~Mengyi Tang who has been the driving force behind Weak IDENT \cite{tang2023weakident} and Fourier IDENT \cite{tang2023weakident}, Prof.~Hongkai Zhao and Prof.~Yimin Zhong for the collaboration on \cite{he2024much}, Dr.~Namjoon Suh, Prof.~Yajun Mei and Prof.~Xiaoming Huo for \cite{he2022asymptotic}, Prof.~Rachel Kuske for \cite{tang2023weakident}, Prof.~Haixia Liu for \cite{he2025group} and Cheng Tang for \cite{tang2025wgidentweakgroupidentification}. Wenjing Liao acknowledges the National Science Foundation under
the NSF DMS 2145167 and the U.S. Department of Energy under the DOE SC0024348. Hao Liu acknowledges National Natural Science Foundation of China grant 12201530 and HKRGC ECS grant 22302123.  

\bibliographystyle{plain}
\bibliography{ref_identR}

\begin{thebibliography}{10}

\bibitem{akaike1974new}
Hirotugu Akaike.
\newblock A new look at the statistical model identification.
\newblock {\em IEEE transactions on automatic control}, 19(6):716--723, 1974.

\bibitem{baake1992fitting}
Ellen Baake, Michael Baake, HG~Bock, and KM~Briggs.
\newblock Fitting ordinary differential equations to chaotic data.
\newblock {\em Physical Review A}, 45(8):5524, 1992.

\bibitem{bellman1969new}
Richard Bellman.
\newblock A new method for the identification of systems.
\newblock {\em Mathematical Biosciences}, 5(1-2):201--204, 1969.

\bibitem{bjorck1990least}
{\AA}ke Bj{\"o}rck.
\newblock Least squares methods.
\newblock {\em Handbook of numerical analysis}, 1:465--652, 1990.

\bibitem{bjorck1991error}
{\AA}ke Bj{\"o}rck.
\newblock Error analysis of least squares algorithms.
\newblock In {\em Numerical Linear Algebra, Digital Signal Processing and
  Parallel Algorithms}, pages 41--73. Springer, 1991.

\bibitem{bongard2007automated}
Josh Bongard and Hod Lipson.
\newblock Automated reverse engineering of nonlinear dynamical systems.
\newblock {\em Proceedings of the National Academy of Sciences},
  104(24):9943--9948, 2007.

\bibitem{brunton2016discovering}
Steven~L Brunton, Joshua~L Proctor, and J~Nathan Kutz.
\newblock Discovering governing equations from data by sparse identification of
  nonlinear dynamical systems.
\newblock {\em Proceedings of the national academy of sciences},
  113(15):3932--3937, 2016.

\bibitem{dai2009subspace}
Wei Dai and Olgica Milenkovic.
\newblock Subspace pursuit for compressive sensing signal reconstruction.
\newblock {\em IEEE transactions on Information Theory}, 55(5):2230--2249,
  2009.

\bibitem{donoho2001uncertainty}
David~L Donoho and Xiaoming Huo.
\newblock Uncertainty principles and ideal atomic decomposition.
\newblock {\em IEEE transactions on Information Theory}, 47(7):2845--2862,
  2001.

\bibitem{fannjiang2012coherence}
Albert Fannjiang and Wenjing Liao.
\newblock Coherence pattern--guided compressive sensing with unresolved grids.
\newblock {\em SIAM Journal on Imaging Sciences}, 5(1):179--202, 2012.

\bibitem{fuchs2004sparse}
J-J Fuchs.
\newblock On sparse representations in arbitrary redundant bases.
\newblock {\em IEEE transactions on Information Theory}, 50(6):1341--1344,
  2004.

\bibitem{gurevich2019robust}
Daniel~R Gurevich, Patrick~AK Reinbold, and Roman~O Grigoriev.
\newblock Robust and optimal sparse regression for nonlinear pde models.
\newblock {\em Chaos: An Interdisciplinary Journal of Nonlinear Science},
  29(10):103113, 2019.

\bibitem{harten1997uniformly}
Ami Harten, Bjorn Engquist, Stanley Osher, and Sukumar~R Chakravarthy.
\newblock Uniformly high order accurate essentially non-oscillatory schemes,
  iii.
\newblock {\em Journal of computational physics}, 131(1):3--47, 1997.

\bibitem{he2025group}
Roy~Y He, Haixia Liu, and Hao Liu.
\newblock Group projected subspace pursuit for block sparse signal
  reconstruction: Convergence analysis and applications.
\newblock {\em Applied and Computational Harmonic Analysis}, 75:101726, 2025.

\bibitem{he2022robust}
Yuchen He, Sung-Ha Kang, Wenjing Liao, Hao Liu, and Yingjie Liu.
\newblock Robust identification of differential equations by numerical
  techniques from a single set of noisy observation.
\newblock {\em SIAM Journal on Scientific Computing}, 44(3):A1145--A1175, 2022.

\bibitem{he2023group}
Yuchen He, Sung~Ha Kang, Wenjing Liao, Hao Liu, and Yingjie Liu.
\newblock Group projected subspace pursuit for identification of variable
  coefficient differential equations (gp-ident).
\newblock {\em Journal of Computational Physics}, 494:112526, 2023.

\bibitem{he2022asymptotic}
Yuchen He, Namjoon Suh, Xiaoming Huo, Sung~Ha Kang, and Yajun Mei.
\newblock Asymptotic theory of-regularized pde identification from a single
  noisy trajectory.
\newblock {\em SIAM/ASA Journal on Uncertainty Quantification},
  10(3):1012--1036, 2022.

\bibitem{he2024much}
Yuchen He, Hongkai Zhao, and Yimin Zhong.
\newblock How much can one learn a partial differential equation from its
  solution?
\newblock {\em Foundations of Computational Mathematics}, 24(5):1595--1641,
  2024.

\bibitem{jaccard1912distribution}
Paul Jaccard.
\newblock The distribution of the flora in the alpine zone. 1.
\newblock {\em New phytologist}, 11(2):37--50, 1912.

\bibitem{kamali2013block}
A~Kamali, MR~Aghabozorgi Sahaf, AM~Doost Hooseini, and AA~Tadaion.
\newblock Block subspace pursuit for block-sparse signal reconstruction.
\newblock {\em Iranian Journal of Science and Technology. Transactions of
  Electrical Engineering}, 37(E1):1, 2013.

\bibitem{kang2021ident}
Sung~Ha Kang, Wenjing Liao, and Yingjie Liu.
\newblock Ident: Identifying differential equations with numerical time
  evolution.
\newblock {\em Journal of Scientific Computing}, 87:1--27, 2021.

\bibitem{kaptanoglu2023benchmarking}
Alan~A Kaptanoglu, Lanyue Zhang, Zachary~G Nicolaou, Urban Fasel, and Steven~L
  Brunton.
\newblock Benchmarking sparse system identification with low-dimensional chaos.
\newblock {\em Nonlinear Dynamics}, 111(14):13143--13164, 2023.

\bibitem{kolmogoroff1936beste}
Andrei Kolmogoroff.
\newblock {\"U}ber die beste ann{\"a}herung von funktionen einer gegebenen
  funktionenklasse.
\newblock {\em Annals of Mathematics}, 37(1):107--110, 1936.

\bibitem{lancaster1981surfaces}
Peter Lancaster and Kes Salkauskas.
\newblock Surfaces generated by moving least squares methods.
\newblock {\em Mathematics of computation}, 37(155):141--158, 1981.

\bibitem{li2020robust}
Jun Li, Gan Sun, Guoshuai Zhao, and H~Lehman Li-wei.
\newblock Robust low-rank discovery of data-driven partial differential
  equations.
\newblock In {\em Proceedings of the AAAI Conference on Artificial
  Intelligence}, volume~34, pages 767--774, 2020.

\bibitem{messenger2021weak}
Daniel~A Messenger and David~M Bortz.
\newblock Weak sindy for partial differential equations.
\newblock {\em Journal of Computational Physics}, 443:110525, 2021.

\bibitem{messenger2021weakode}
Daniel~A Messenger and David~M Bortz.
\newblock Weak sindy: Galerkin-based data-driven model selection.
\newblock {\em Multiscale Modeling \& Simulation}, 19(3):1474--1497, 2021.

\bibitem{muller2004parameter}
Thomas~G M{\"u}ller and Jens Timmer.
\newblock Parameter identification techniques for partial differential
  equations.
\newblock {\em International Journal of Bifurcation and Chaos},
  14(06):2053--2060, 2004.

\bibitem{reinbold2020using}
Patrick~AK Reinbold, Daniel~R Gurevich, and Roman~O Grigoriev.
\newblock Using noisy or incomplete data to discover models of spatiotemporal
  dynamics.
\newblock {\em Physical Review E}, 101(1):010203, 2020.

\bibitem{rudy2019data}
Samuel Rudy, Alessandro Alla, Steven~L Brunton, and J~Nathan Kutz.
\newblock Data-driven identification of parametric partial differential
  equations.
\newblock {\em SIAM Journal on Applied Dynamical Systems}, 18(2):643--660,
  2019.

\bibitem{Rudy2019DataDrivenIO}
Samuel~H. Rudy, Alessandro Alla, Steven~L. Brunton, and J.~Nathan Kutz.
\newblock Data-driven identification of parametric partial differential
  equations.
\newblock {\em SIAM J. Appl. Dyn. Syst.}, 18:643--660, 2019.

\bibitem{rudy2017data}
Samuel~H Rudy, Steven~L Brunton, Joshua~L Proctor, and J~Nathan Kutz.
\newblock Data-driven discovery of partial differential equations.
\newblock {\em Science Advances}, 3(4):e1602614, 2017.

\bibitem{schaeffer2017learning}
Hayden Schaeffer.
\newblock Learning partial differential equations via data discovery and sparse
  optimization.
\newblock {\em Proceedings of the Royal Society A: Mathematical, Physical and
  Engineering Sciences}, 473(2197):20160446, 2017.

\bibitem{schmidt2009distilling}
Michael Schmidt and Hod Lipson.
\newblock Distilling free-form natural laws from experimental data.
\newblock {\em science}, 324(5923):81--85, 2009.

\bibitem{stone1978cross}
Mervyn Stone.
\newblock Cross-validation: A review.
\newblock {\em Statistics: A Journal of Theoretical and Applied Statistics},
  9(1):127--139, 1978.

\bibitem{tang2025wgidentweakgroupidentification}
Cheng Tang, Roy~Y. He, and Hao Liu.
\newblock Wg-ident: Weak group identification of pdes with varying
  coefficients, 2025.

\bibitem{tang2023weakident}
Mengyi Tang, Wenjing Liao, Rachel Kuske, and Sung~Ha Kang.
\newblock Weakident: Weak formulation for identifying differential equation
  using narrow-fit and trimming.
\newblock {\em Journal of Computational Physics}, 483:112069, 2023.

\bibitem{tang2023fourier}
Mengyi Tang, Hao Liu, Wenjing Liao, and Sung~Ha Kang.
\newblock Fourier features for identifying differential equations
  (fourierident).
\newblock {\em arXiv preprint arXiv:2311.16608}, 2023.

\bibitem{tibshirani1996regression}
Robert Tibshirani.
\newblock Regression shrinkage and selection via the lasso.
\newblock {\em Journal of the Royal Statistical Society Series B: Statistical
  Methodology}, 58(1):267--288, 1996.

\bibitem{tropp2004just}
Joel~A Tropp.
\newblock Just relax: Convex programming methods for subset selection and
  sparse approximation.
\newblock {\em ICES report}, 404, 2004.

\bibitem{tropp2006just}
Joel~A Tropp.
\newblock Just relax: Convex programming methods for identifying sparse signals
  in noise.
\newblock {\em IEEE Transactions on Information Theory}, 52(3):1030--1051,
  2006.

\bibitem{yuan2006model}
Ming Yuan and Yi~Lin.
\newblock Model selection and estimation in regression with grouped variables.
\newblock {\em Journal of the Royal Statistical Society: Series B (Statistical
  Methodology)}, 68(1):49--67, 2006.

\end{thebibliography}

\end{document}